\documentclass{article}

\usepackage{amssymb}
\usepackage{latexsym}

\def \IQ {\mathbb Q}
\def \IR {\mathbb R}
\def \IN {\mathbb N}

\def \proof{\noindent {\bf Proof\ \  }}

\def \range {{\rm range}}

\def \mod {{\rm mod}}
\def \an {\mbox{\ \ and\ \ }}

\def \In{\subseteq}

\newcommand{\qq}{\phantom{.}\hfill$\Box$}

\newcommand{\mto}{\rightrightarrows}
\newcommand{\mmto}{\mbox{
\setlength{\unitlength}{1em}
\begin{picture}(0.4,0)
\makebox(0,0.6){$\mbox{\scriptsize \raisebox{0.083em}{$|$}}
\hspace*{-1.1ex}\mto$}
\end{picture}
}}

\newtheorem{definition}{Definition}[section]{\bf}{\it}
\newtheorem{theorem}[definition]{Theorem}{\bf}{\it}
{\bf}{\it}
\newtheorem{lemma}[definition]{Lemma}{\bf}{\it}
\newtheorem{corollary}[definition]{Corollary}{\bf}{\it}
{\bf}{\it}
\newtheorem{example}[definition]{Example}{\bf}{\it}
{\bf}{\it}
{\bf}{\it}

\begin{document}

\author{Klaus Weihrauch}
\title{Computable paths intersect in a computable point}

\maketitle

\begin{abstract} Consider two paths $f,g:[0;1]\to [0;1]^2$ in the unit square such that $f(0)=(0,0)$, $f(1)=(1,1)$, $g(0)=(0,1)$ and $g(1)=(1,0)$. By continuity of $f$ and $g$ there is a point of intersection. We prove that
there is a computable point of intersection if $f$ and $g$ are computable.
\end{abstract}

\section{Introduction}\label{seca}

%\frame{
%\vspace*{30ex}\hspace*{20ex}
%$\includegraphics[width=.3\linewidth]{montreal-turing.jpg}$
%}

%\begin{center}
%$\includegraphics[width=.6\linewidth]{montreal-turing.jpg}$
%\end{center}

A {\em path} in the Euclidean plane is a continuous function $f: [0;1]\to \IR^2$, a {\em curve} is the range of a path. The following is known about planar curves.

\begin{theorem}\label{t1}
Let  $f,g:[0;1]\to [0;1]^2$  be two paths in the unit square such that
\begin{eqnarray}
\label{f1}& f(0)=(0,0), \  f(1)=(1,1), \   g(0)=(0,1)\  \mbox{ and } \ g(1)=(1,0) \,.
\end{eqnarray}
Then the two curves $\range(f)$ and $\range(g)$ intersect.
\end{theorem}

\begin{figure}[htbp]
\setlength{\unitlength}{1.0pt}
\linethickness{0.7pt}
\begin{picture}(200,140)(-225,-20)

\put(0,0){\vector (1,0){110}}
\put(0,100){\line (1,0){100}}
\put(0,0){\vector (0,1){110}}
\put(100,0){\line (0,1){100}}
\qbezier(0,0)(20,80)(100,100)
\qbezier(0,100)(50,10)(100,0)

\put(54,70){\makebox(0,0)[cc]{$f$}}
\put(40,30){\makebox(0,0)[cc]{$g$}}

\put(-28,-8){\makebox(0,0)[cc]{$f(0)=(0,0)$}}
\put(-28,108){\makebox(0,0)[cc]{$g(0)=(0,1)$}}
\put(128,-8){\makebox(0,0)[cc]{$(1,0)=g(1)$}}
\put(128,108){\makebox(0,0)[cc]{$(1,1)=f(1)$}}

\end{picture}
\caption{Intersecting curves} \label{fig1}
\end{figure}
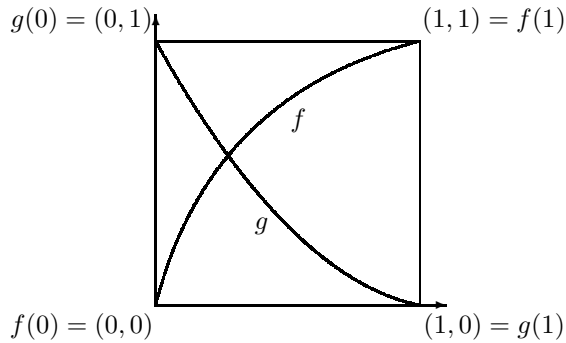

Figure~\ref{fig1} visualizes the theorem.
Manukyan \cite{Man76a} (in Russian) has proved the following surprising theorem, cited in \cite[Page 279]{Kus99} as follows:
\begin{theorem}[Manukyan]\label{t2}There are two constructive (and therefore continuous) planar curves $ \varphi_1$ and $ \varphi_2$ such that
\begin{eqnarray*}
& \varphi_1(0)=(0,0), \  \varphi_1(1)=(1,1),\
 \varphi_2(0)=(0,1), \   \varphi_2(1)=(1,0)\,, \\
&\mbox{for every $0<t< 1$ both $ \varphi_1(t)$  and $ \varphi_2(t)$ belong to the open unit square,}\\
&\mbox{ $ \varphi_1$ and $ \varphi_2$ do not intersect.}
\end{eqnarray*}
\end{theorem}

Since the paths in Theorem~\ref{t1} are continuous, the two theorems seem to be contradictory. The contradiction vanishes by the observation that
 in Theorem~\ref{t2}, ``constructive'' means  ``Markov computable'' or computable in the  ``Russian'' way   \cite{Kus73b,Abe80,Wei00}. In the Russian approach only the countable set $\IR_c$ of computable real numbers is considered, hence $ \varphi_1,  \varphi_2:\IR_c\cap[0;1]\to (\IR_c\cap[0;1])^2$.
\smallskip

In this article we use the definition of computable real functions $f$ and $g$ introduced by Grzegorczyk and Lacombe \cite{Grz55,Grz57,Lac55a} and call them {\em computable}. We assume that the reader is familiar with the basic properties of this concept \cite{Wei00,BHW08,WG09}.
We will prove the following theorem.
\begin{theorem}\label{t3}
Let  $f,g:[0;1]\to [0;1]^2$  be two {\bf computable} paths in the unit square such that
\begin{eqnarray}
\label{f3}& f(0)=(0,0), \  f(1)=(1,1), \   g(0)=(0,1)\  \mbox{ and } \ g(1)=(1,0) \,.
\end{eqnarray}
Then there is a computable point $x\in \range(f)\cap \range(g)$.
\end{theorem}
The restrictions  $\varphi_f$ and $\varphi_g$ of the computable functions $f$ and $g$, respectively, to the computable points are Markov-computable \cite{Wei00}. This is is no contradiction to Theorem~\ref{t2} but merely means that the functions $\varphi_1$ and $\varphi_2$ from Theorem~\ref{t2} have no extensions to computable functions $[0;1]\to[0;1]^2$.
The existence of computable intersection points in various other situations has been studied by   Iljazovi\'{c} and  and Pa\v{z}ek in \cite{IP17}.

In the next sections we will prove the following formally weaker but equivalent version of Theorem~\ref{t3}.
\begin{theorem}\label{t6}
Let  $f,g:[0;1]\to [0;1]^2$  be two {\bf computable} paths in the unit square such that
\begin{eqnarray}
\label{f34}& f(0)=(0,0), \  f(1)=(1,1), \   g(0)=(0,1), \ g(1)=(1,0) \ \mbox{and}\\
 \label{f35}&\mbox{for every $0<t< 1$ both $f(t)$  and $g(t)$ belong to the open unit square.}
\end{eqnarray}
Then there is a computable point $x\in \range(f)\cap \range(g)$.
\end{theorem}
Obviously Theorem~\ref{t3} implies  Theorem~\ref{t6}. The converse is also true: Suppose, Theorem~\ref{t6} is true. Let $f,g:[0;1]\to [0;1]^2$ be computable functions such that (\ref{f3}) is true.
We extend $f$ and $g$ to functions $f_1,g_1:[-1;2]\to [-1;2]^2$ by
$f_1(t):=(t,t)$ and $g_1(t):=(t,1-t)$ for $-1\leq t <0$ and $1<t\leq 2$.
We transform the square $[-1;2]^2$ together with $f_1$ and $ g_1$ to the unit square. For $0\leq t\leq 1$ let
$f_2(t):=(f_1(3t-1)+(1,1))/3$ and $g_2(t):=(g_1(3t-1)+(1,1))/3$.
Then $f_2$ and $g_2$ satisfy  (\ref{f34}) and (\ref{f35}).
By Theorem~\ref{t6} there are (not necessarily computable) numbers $s,t\in[0;1]$ such that $f_2(s)=x= g_2(t)$ and $x\in \IR^2$ is computable, hence $f_1(3s-1)=3x-(1,1)=g_1(3t-1)$.
Since $f_1$ and $g_1$ cannot intersect outside  the unit square
$f(3s-1)=3x-(1,1)=g(3t-1)$ where $0\leq 3s-1\leq 1$ and $0\leq 3t-1\leq 1$.
Therefore $3x-(1,1)\in\range(f)\cap \range(g)$. The point $3x-(1,1)$ is computable.

A corollary of Theorem~\ref{t6} is the computable intermediate value theorem~\cite{Wei00}.
\begin{theorem}[computable intermediate value]\label{t8}
 Every computable function $h:[0;1]\to [-1;1]$ with $h(0)=-1$ and $h(1)=1$ has a computable zero (that is,  $h(t)=0$ for some computable number $t$). \end{theorem}
For a proof as a corollary of Theorem~\ref{t6}, define $f(t):=(t,0)$ and $g(t):=(t,h(t))$. By an easily provable corollary of Theorem~\ref{t6}, there are numbers $s,t$ and a computable point $(x,y)$ such that $f(s)=(x,y)=g(t)$, hence
$(s,0)=(x,y)=(t,h(t))$. Therefore, $t=x$ is computable and $h(t)=0$.

The intermediate value theorem can be proved directly as follows:
\begin{enumerate}
\item[{\bf A:}] \label{en1}
Define a predicate $Q_{iv}$ by
\begin{eqnarray}\label{f61}Q_{iv}\iff h(a;b)=\{0\}\ \ \mbox{for some rational numbers} \ \ a<b\,.
\end{eqnarray}
\item[{\bf B:}] \label{en2}
 Suppose $Q_{iv}$. Then $f(c)=0$ for some rational (hence computable) number $c\in (a;b)$.
\item[{\bf C:}] \label{en3}
Suppose $\neg Q_{iv}$. Call an interval $I=(a;b)$ a crossing, if $a,b\in\IQ$ and $h(a)<0<h(b)$. By continuity of $h$ the following can be shown easily:
\begin{eqnarray}\label{f62}\mbox{For every crossing $I$ there is a crossing $J$ such that $\overline J\In I$ and ${\rm length}(J)< {\rm length}(I)/2$\,.}
\end{eqnarray}
\item[{\bf D:}] \label{en4} Since the set of crossings is c.e., beginning with $(0;1)$  we can compute a nested sequence of crossings which converges to some real number $t$. Then $h(t)=0$.
\end{enumerate}
We will prove Theorem~\ref{t6} in similar steps.

Computable curves can be extremely complicated.  Consider, for example, a space-filling curve or a curve with infinitely many spirals, each of which containing infinitely many sub-spirals etc. infinitely often or curves with chaotic behavior. Furthermore, there is a computable function $h:[0;1]\to [-1;1]$ with $h(0)=-1$ and $h(1)=1$ such that $h^{-1}\{0\}$ has measure $>1/2$ but contains only a single computable number.\footnote{
There is a computable open set $U\In \IR$ with measure $<1/2$ which contains the computable real numbers \cite{Spe59}\cite[Theorem~4.2.8]{Wei00}. Define a computable function $h_0:\IR\to\IR$ such that $(\forall t)h_0(t)\leq 0$ and $h_0^{-1}(0)=\IR\setminus U$. Construct $h$ from $h_0$.}
Therefore at first glance there might be enough freedom to construct computable functions $f$ and $g$ which avoid to cross in a computable point.

In Sections~\ref{secm} to~\ref{sece} we prove  Theorem~\ref{t6} corresponding to Steps {\bf A,B,C,D} above.
In Section~\ref{secm} we define a predicate $Q$ which corresponds to the
$Q_{iv}$ in Step~A above.
In Section~\ref{secb} we prove that $Q$ implies the existence of a computable point of intersection. In
Section~\ref{secc} assuming $\neg Q$ we define crossings and show that every crossing has a much smaller sub-crossing.
In section~\ref{sece} we prove that every crossing has a proper sub-crossing and that the set of proper crossings is c.e. Then we compute a point of intersection of $f$ and $g$ as the limit of a decreasing sequence of proper crossings.
In Section~\ref{secf} we outline a proof for the special case that $f$ is injective and add some further considerations.

Computability on the real numbers, on open  and on compact subsets $\IR$ and of $\IR^2$ is defined via canonical
representations\footnote{A name of an open set $U$ is a list $B_1,B_2,\ldots $ of open rational balls such that
$U=\bigcap B_i$, a name of a compact set $K$ is a list of all finite covers of $K$ with rational balls.} \cite{Wei00,WG09}.
In the following, $f$ and $g$ will be computable functions satisfying the conditions (\ref{f34}) and (\ref{f35}) from Theorem~\ref{t6}.

\section{Two alternatives}\label{secm}

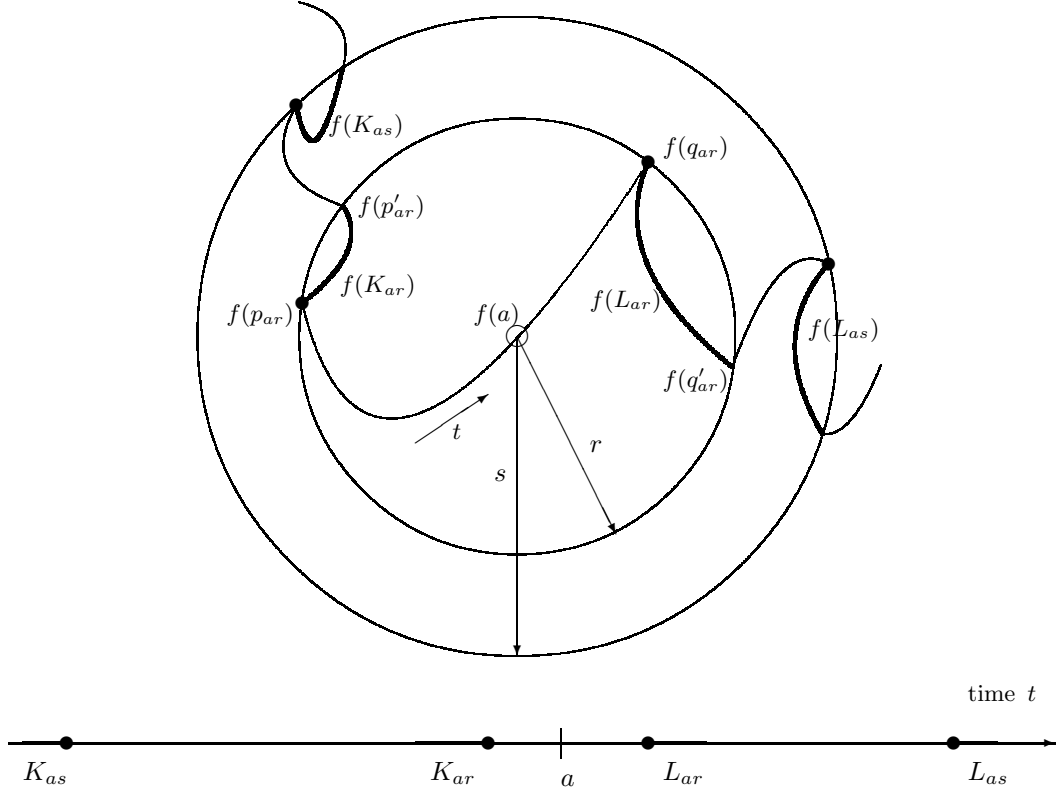
\begin{figure}[htb]
\setlength{\unitlength}{5.5pt}
\begin{picture}(80,54)(-12,-12)

\newsavebox{\ballx}
\savebox{\ballx}{
\qbezier(15,0)(15,6.21)(10.61,10.61)
\qbezier(10.61,10.61)(6.21,15)(0,15)
\qbezier(15,0)(15,-6.21)(10.61,-10.61)
\qbezier(10.61,-10.61)(6.21,-15)(0,-15)
\qbezier(-15,0)(-15,6.21)(-10.61,10.61)
\qbezier(-10.61,10.61)(-6.21,15)(0,15)
\qbezier(-15,0)(-15,-6.21)(-10.61,-10.61)
\qbezier(-10.61,-10.61)(-6.21,-15)(0,-15)
}

\put(35,20){\usebox{\ballx}\circle{1.5}}
\put(35,20){\vector(1,-2){6.8}}
\put(40,12){\parbox{10ex}{ $r$ }}

\newsavebox{\ballw}
\savebox{\ballw}{
\qbezier(22,0)(22,9.11)(15.56,15.56)
\qbezier(15.56,15.56)(9.11,22)(0,22)
\qbezier(22,0)(22,-9.11)(15.56,-15.56)
\qbezier(15.56,-15.56)(9.11,-22)(0,-22)
\qbezier(-22,0)(-22,9.11)(-15.56,15.56)
\qbezier(-15.56,15.56)(-9.11,22)(0,22)
\qbezier(-22,0)(-22,-9.11)(-15.56,-15.56)
\qbezier(-15.56,-15.56)(-9.11,-22)(0,-22)
%\put(0,0){\circle*{.5}}
}
\put(35,20){\usebox{\ballw}}
\put(35,20){\vector(0,-3){22}}
\put(33.4,10){\parbox{10ex}{ $s$ }}

\put(32,21){\parbox{1ex}{\small $f(a)$ }}

\qbezier(20.2,22.3)(25,2.5)(44,32)
%\qbezier(44,32)(41,25)(49.8,17.9)
%\qbezier(20.2,22.3)(25,26)(23,29)
\qbezier(23,29)(17,31)(19.8,35.9)
%\qbezier(19.8,35,9)(21,30)(23,38.47)
\qbezier(23,38.47)(24,42)(20,43)
\qbezier(49.8,17.9)(53,27)(56.4,25)
%\qbezier(56.4,25)(52,20)(56,13.3)
\qbezier(56,13.3)(58,13)(60,18)

\linethickness{1.3pt}
\qbezier(19.8,35,9)(21,30)(23,38.47)
\qbezier(44,32)(41,25)(49.8,17.9)
\qbezier(20.2,22.3)(25,26)(23,29)
\qbezier(56.4,25)(52,20)(56,13.3)

\thinlines
\put(0,-8){\vector(1,0){72}}
\put(38,-7){\line(0,-1){2}}

\put(66,-5.){\parbox{6ex}{\small time \ $t$ }}

\linethickness{1.3pt}
\put(1,-8){\line(1,0){3}}
\put(28,-8){\line(1,0){5}}
\put(44,-8){\line(1,0){4}}
\put(65,-8){\line(1,0){3}}
\thinlines

\put(1,-10.5){\parbox{10ex}{ $K_{as}$ }}
\put(29,-10.5){\parbox{10ex}{ $K_{ar}$ }}
\put(38,-11.){\parbox{10ex}{ $a$ }}
\put(45,-10.5){\parbox{10ex}{ $L_{ar}$ }}
\put(66,-10.5){\parbox{10ex}{ $L_{as}$ }}

\put(20.2,22.3){\circle*{.9}}
\put(44,32){\circle*{.9}}
\put(19.8,35,9){\circle*{.9}}
\put(56.4,25){\circle*{.9}}

\put(4,-8){\circle*{.9}}
\put(33,-8){\circle*{.9}}
\put(44,-8){\circle*{.9}}
\put(65,-8){\circle*{.9}}

\put(28,12.7){\vector(3,2){5}}
\put(30,13){\parbox{2ex}{\small \ $t$ }}

\put(15,21){\parbox{1ex}{\small $f(p_{ar})$ }}
\put(45,32.5){\parbox{1ex}{\small $f(q_{ar})$ }}
\put(24,28.5){\parbox{1ex}{\small $f(p'_{ar})$ }}
\put(45,16.5){\parbox{1ex}{\small $f(q'_{ar})$ }}
\put(22.8,23){\parbox{1ex}{\small $f(K_{ar})$ }}
\put(40,22){\parbox{1ex}{\small $f(L_{ar})$ }}

\put(22.,34){\parbox{1ex}{\small $f(K_{as})$ }}
\put(55,20){\parbox{1ex}{\small $f(L_{as})$ }}

\end{picture}
\caption{A segment of ${\rm graph}(f)$ around time $a$ crossing balls $B(f(a),r)$ and $B(f(a),s)$. }\label{fig5}
\end {figure}

In this section we define a predicate $Q$ which corresponds to the predicate $Q_{iv}$ from Step~{\bf A} in Section~\ref{seca}.
We consider Euclidean balls $B(f(a),r))$ such that $a,r\in\IQ$ and
$\overline{B(f(a),r))}\In(0;1)^2$.
Since $\overline{B(f(a),r))}\In(0;1)^2$ iff $0<a<1$, \ $r>0$ and ${\rm pr}_i\circ f(a)\in (r;1-r)$ (for $i=1,2$),
\begin{eqnarray}\label{f4}
\{(a,r)\in\IQ^2\mid \overline{B(f(a),r))}\In(0;1)^2\} & \mbox{is c.e.}
\end{eqnarray}
For $a,r\in\IQ$ such that $\overline{B(f(a),r))}\In(0;1)^2$ define
\begin{eqnarray}
\label{f11} (p_{ar};q_{ar})& :=&\mbox{the longest open interval $I$ such that }\ a\in I \ \mbox{ and }\ f(I)\In B(f(a),r)\,,\   \\
\label{f13} (p_{ar}';q_{ar}')&:=&  \mbox{the longest open interval $I$ such that }\ a\in I \ \mbox{ and }\ f(\overline I)\In \overline{B(f(a),r)}\,, \\
\label{f7}K_{ar} & := & [p'_{ar};p_{ar}]\,,\\
\label{f15}L_{ar} & := & [q_{ar};q'_{ar}]\,.
\end{eqnarray}
The definitions are illustrated in Figure~\ref{fig5}.
%
%
%\begin{figure}[htbp]
We can interpret the arguments of $f$ and $g$ as ``time''.
Starting at time $a$ in positive direction $f$ leaves the open ball at time $q_{ar}$ and leaves the closed ball at time $q'_{ar}$. Correspondingly, starting at time $a$ in negative direction $f$ leaves the open ball at time $p_{ar}$ and leaves the closed ball at time $p'_{ar}$.
For radii $r$ and $s$, Figure~\ref{fig5} shows a time interval from $t<\min K_{as}$ to $t>L_{as}$ and its image.
 The function $f$ can be much more complicated than shown in the figure. For example there is a computable function $f$ such that $f(p_{ar};q_{ar})=B(f(a),r)$ and $f(K_{ar})=f(L_{ar})=\overline{B(f(a),r)}$.
Notice that $f$ may cross the ball $B(f(a),r)$ additionally at times $t$ with $t<p_{ar}'$ or $t>q_{ar}'$.

Obviously, $f(p_{ar}),f(q_{ar}),f(p'_{ar}),f(q'_{ar})\in
 \partial(B(f(a),r)$
\footnote{where $\partial V$ denotes the boundary of $V$} and
\begin{eqnarray}
\label{f13a} p_{ar} &= & \inf\{ b\in (0;a)\cap \IQ  \mid f[b;a]\In B(f(a),r)\}\,, \\
\label{f11a} q_{ar} &= & \sup\{b\in (a;1)\cap \IQ\mid f[a;b]\In B(f(a),r)\}\,,\\
\label{f14a} p_{ar}'&= & \sup\{b\in (0;a)\cap \IQ\mid f(b)\not\in \overline {B(f(a),r)}\}\,,\\
\label{f12a} q_{ar}' &= & \inf\{b\in (a;1)\cap \IQ\mid f(b)\not\in \overline {B(f(a),r)}\}\,.
\end{eqnarray}
The numbers $q_{ar},q_{ar}',p_{ar}$ and $p_{ar}'$ are  not computable in general but semi-computable.
\begin{lemma}\label{l4}
The sets
\vspace*{-1ex}
\begin{eqnarray}
\label{f5}\{(a,r,c,d)\in\IQ^4\mid \overline{B(f(a),r))}\In(0;1)^2
\ \mbox{ and } \ K_{ar}\In (c;d)\} & \mbox{and}\\
\label{f6}\{(a,r,c.d)\in\IQ^4\mid \overline{B(f(a),r))}\In(0;1)^2
\ \mbox{ and } \ L_{ar}\In (c;d)\}
\end{eqnarray}
\vspace*{-1ex}
are c.e.\\[-2ex]
$ $
\end{lemma}

\proof
From the compact set $[a;b]$ we can compute the compact set $f[a;b]$. Since $K\In O$ for compact $L$ and open $O$ is c.e.  we can  enumerate the set of all $b\in (a;1)\cap \IQ$ such that $f[a;b]\In B(f(a),r)$.
Since $ f(b)\not\in \overline {B(f(a),r)}$ is equivalent to
$\|f(a)-f(b)\|<r$ (in $\IR^2$), we can enumerate the set of all $b\in (a;1)\cap\IQ$ such that $ f(b)\not\in \overline {B(f(a),r)}$.
(\ref{f6}) follows from (\ref{f15}),( \ref{f11a}) and (\ref{f12a}).
Accordingly $K_{ar}$.
\qq\\

We define a predicate $Q$ which in this proof corresponds to the predicate $Q_{iv}$ from the proof of the computable intermediate value theorem in Section~\ref{seca}.
\begin{eqnarray}\label{f43}
Q:\Leftrightarrow\left\{ \begin{array}{l}
\mbox{ there are $ a,r,s\in\IQ)$ with $0<r<s$ and   $\overline{B(f(a),s)}\In (0;1)^2)$ such that} \\
(\forall t\in\IQ\cap[r;s])\ ( f\circ\max(K_{at})\in \range(g) \ \mbox{ or } \ f\circ\min(L_{at})\in \range(g))\,.
\end{array}\right\}
\end{eqnarray}
In Figure~\ref{fig5}, $ f\circ\max(K_{ar})$ and $f\circ\min(L_{ar})$ are the marked points on the circle with radius $r$. By $Q$ for every radius
$r\leq t \leq s$ at least one of the corresponding marked points must be in $\range(g)$.
We will prove that $f(s)=f(t)$ for some computable $s$ and $t$  if $Q$ and if $\neg Q$. In general neither $Q$ nor $\neg Q$ are c.e.

\section{The first alternative}\label{secb}
We assume $Q$  from (\ref{f43}). Therefore,  there are $ a,r,s\in\IQ)$ with $0<r<s$ and   $\overline{B(f(a),s)}\In (0;1)^2)$ such that
\begin{eqnarray}\label{f9}
&(\forall t\in\IQ\cap[r;s])\ ( f\circ\max(K_{at})\in \range(g) \ \mbox{ or } \ f\circ\min(L_{at})\in \range(g))\,.
\end{eqnarray}

For an interval $I$ let $lg(I)$ be its length.
 For closed real intervals $K,L$ define $K<L$ iff $\max(K)<\min(L)$. If $K<L$ let $lg(K\cup L):= \max(L)-\min(K)$ \footnote{We consider $(K,L)$ as a generalized interval with lower bound $K$ and upper bound $L$.}. Then (see Figure~\ref{fig5})
\begin{eqnarray}\label{f2} K_{as}<K_{ar}<L_{ar}<L_{as}\ \ & \ \ \mbox{if} \ \ & r<s\,.
\end{eqnarray}
\begin{lemma}\label{l5} From rational $a,r,s,\varepsilon $ such that $0<r<s$, \ $\overline{B(f(a),s)}\In (0;1)^2$
and $\varepsilon >0$ we can compute rational numbers \  $\overline r,\overline s$ and rational intervals $I,J$ such that
\begin{eqnarray}\label{f12}
 r<\overline r<\overline s<s,  \ \ \
 K_{a\overline s}\cup K_{a\overline r}\In I,\ \ \
 L_{a\overline r}\cup L_{a\overline s}\In J,\ \ \
 lg(I)<\varepsilon \ \ \mbox{and}\ \ \ lg(J)<\varepsilon
\,.
\end{eqnarray}
\end{lemma}

\proof First we prove that  such numbers $\overline r,\overline s$ exist.
There is a sequence $(r_i)_{i\in\IN}$ of rational numbers such that
$r<r_0<r_1<r_2<\ldots <s$. Then by (\ref{f2})
$$K_{as}<\ldots<K_{ar_2}<K_{ar_1}<K_{ar_0}<K_{ar} \ \ \  < \ \ \
L_{ar}<L_{ar_0}<L_{ar_1}<L_{ar_2}<\ldots< L_{as}\,.$$
Since all these intervals are pairwise disjoint
There are pairwise disjoint rational intervals $I_i$ and $J_i$ such that
$K_{ar_{2i+1}}\cup K_{ar_{2i}}\In I_i$ and
$L_{ar_{2i}}\cup L_{ar_{2i+1}}\In J_i$.
Since $\sum_ilg(I_i) + \sum_i lg(J_i) \leq \max(L_{as})-\min(K_{as})$
there are only finitely many numbers $i$ such that $lg(I_i)\geq \varepsilon$ and only finitely many numbers $i$ such that $lg(J_i)\geq \varepsilon$. Hence for some number $i$, $lg(I_i)<\varepsilon$ and
$lg(J_i)<\varepsilon$.
Choose $\overline r:=r_{2i}$, $\overline s:=r_{2i+1}$, $I:=I_{2i}$ and $J:=J_{2i}$.

In general, for closed intervals $K<L$,
$K\cup L\In I$  iff there are disjoint rational intervals $I_1,I_2$ such that $K\In I_1$, $L\In I_2$ and $I_1\cup I_2\In I$.
Therefore by Lemma~\ref{l4} we can compute appropriate radii $\overline r,\overline s$ and intervals $I$ and $J$.
\qq

\begin{theorem}\label{t4} If $Q$ then $f(\alpha)\in \range(g)$ for some computable number $\alpha\in\IR$.
\end{theorem}

\proof
By Lemma~\ref{l5} from $a,r,s$ we can compute a sequence $(r_n,s_n)_n$ of pairs of rational radii and sequences $(I_n)_n$ and $(J_n)_n$ of rational intervals such that $r<r_0<r_1<r_2<\ldots <s_2<s_1<s_0<s$ and for all $n$,
\begin{eqnarray}\label{f14}
\lg(I_n)<2^{-n}\ ,\ \ \  \lg(J_n)<2^{-n}\,,\ \ \
K_{as_n}\cup K_{ar_n}\In I_n\,,\ \ \  L_{ar_n}\cup L_{as_n}\In J_n\,.
\end{eqnarray}

Let
$$S_n:=\bigcap_{m=0}^n \overline I_m \ \ \ \mbox{and}\ \ \ T_n:=\bigcap_{m=0}^n \overline J_m$$

Suppose $m<n$. Then $r_m<r_n<s_n<s_m$, hence $L_{ar_m}<L_{ar_n}<L_{as_n}<L_{as_m}$, therefore, $L_{ar_n}\cup L_{as_n}\In J_m$.

Therefore, $(T_n)_n$ is a descending sequence of closed rational intervals such that $L_{ar_n}\cup L_{as_n}\In T_n$ and $lg(T_n)<2^{-n}$ by (\ref{f14}). The sequence converges to a point $\beta\in\IR$ which is computable since the sequence $(J_n)_n$ is computable.
Correspondingly, $K_{as_n}\cup K_{ar_n}\In S_n$, $\lg(S_n)<2^{-n}$ and the sequence $(S_n)_n$ converges to a computable point $\alpha$.

We apply assumption $Q$. By (\ref{f9}),
\begin{eqnarray}\label{f59}
&(\forall n)\ ( f\circ\max(K_{ar_n})\in \range(g) \ \mbox{ or } \ f\circ\min(L_{ar_n})\in \range(g))
\end{eqnarray}
(see the marked points in Figure~\ref{fig5}), hence
\begin{eqnarray}\label{f60}f\circ\max(K_{ar_n})\in \range(g) \ \mbox{infinitely often}\ \ \ \ \ \mbox{or}\ \ \ \ \ f\circ\min(L_{ar_n})\in \range(g)\ \  \mbox{infinitely often}\,.
\end{eqnarray}

Suppose $f\circ\max(K_{ar_n})\in \range(g)$ infinitely often.
Then for some increasing sequence $n_0, n_1,n_2,\ldots$ of indices, $f(\alpha_i)\in\range(g)$ where $\alpha_i:= \max(K_{ar_{n_i}})$. By $i\leq n_i$, $\alpha_i\in K_{ar_{n_i}}\In S_{n_i}\In S_i$,
hence the sequence $(\alpha_i)_i$ converges to $\alpha$.
Since $f$ is continuous, the sequence $(f(\alpha_i)_i$ converges to $f(\alpha)$. Since $f(\alpha_i)\in\range(g)$, $(f(\alpha_i)_i$ is a sequence in $\range(g)$ which converges to $f(\alpha)$. Since $\range(g)$ is compact it is complete. Therefore, $f(\alpha)\in \range(g)$.
Correspondingly, if $f\circ\min(L_{ar_n})\in \range(g)$ infinitely often then $f(\beta)\in\range(g)$.

Therefore, we have shown that $Q$ implies $f(\alpha)\in\range(g)$ for some computable number $\alpha\in\IR$.
\qq

\begin{corollary}\label{cor1}
If $Q$ then $g(\alpha)\in\range(f)$ for some computable number $\alpha\in\IR$.
\end{corollary}
\proof By symmetry.
\qq

\section{Gates}\label{secc}
In this section we generalize Step~{\bf C} of the proof of the computable intermediate value theorem in Section~\ref{seca}. We define gates and prove that every gate crossed by $g$  contains a much smaller gate crossed by $g$. In Section~\ref{sece} we compute an intersection point of $f$ and $g$ as the limit of a converging sequence of crossings. Let $Q$ be the predicate defined in (\ref{f43}). In the following we assume $\neg Q$, that is,
for all $a,r,s\in\IQ$ such that $0<r<s$ and   $\overline{B(f(a),s)}\In (0;1)^2)$,
\begin{eqnarray}\label{f23}(\exists\ t\in\IQ\cap(r;s))\ ( f\circ\max(K_{at})\not\in \range(g) \ \mbox{ and } \ f\circ\min(L_{at})\not\in \range(g))\,.
\end{eqnarray}

First we introduce gates and a concept of crossing. Gates generalize the intervals in our proof of the computable intermediate value theorem in Section~\ref{seca}.
\begin{definition}[passage, gate]\label{d3}$ $
\vspace*{-1.5ex}
\begin{enumerate}
\item \label{d3d}A {\rm passage} is a pair $(V,c)$  such that
   $V$ is is the open unit square $(0;1)^2$ or a ball $B(f(a),r)$ ($a,r\in\IQ$) with $\overline{B(f(a),r)}\In (0;1)^2$ or the
   strict intersection\footnote
   {where we call the intersection of $B_1$ and $B_2$ strict if $B_1\not\In B_2$, \ $B_2\not\In B_1$ and $B_1\cap B_2\neq\emptyset$}
   of two such balls, and
   $c\in\IQ\cap(0;1)$ with $f(c)\in V$.

\item \label{d3a} Let  $I_{Vc}$ be the longest open interval $I$ such that $c\in I$ and $f(I)\In V$ and for $d\in\IQ$ with $g(d)\in V$ let $I^g_{Vd}$ be the longest open interval $I$ such that $d\in I$ and $g(I)\In V$ (cf. (\ref{f11}) and [\ref{f13})).

\item \label{d3e} A gate is a passage such that
$f\circ \inf (I_{Vc})\not\in\range(g)$ and $f\circ sup (I_{Vc})\not\in\range(g)$.

\item \label{d3b}
 A gate $(V,c)$ is crossed by $g$ via $d\in\IQ$ if
   $g(d)\in V$ and the straight line segments
from   $f\circ \inf (I_{Vc})$ to  $f\circ \sup (I_{Vc})$ and
 from   $f\circ \inf (I^g_{Vd})$ to  $f\circ \sup (I^g_{Vd})$
 intersect. We say that the gate $(V,c)$ is crossed by $g$ if it is crossed via some  $d$.
 %A gate is a crossing, if it is crossed by $g$ via some $d$.
\end{enumerate}
\end{definition}
Notice that by (\ref{f7}) and (\ref{f15}),
\begin{eqnarray}\label{f20}
\inf(I_{Va})=\max(K_{ar}) & \mbox{and} &
\sup(I_{Va})=\min(L_{ar}) \ \ \ \mbox{if}\ \ \  V=(B(f(a),r),a)\,.
\end{eqnarray}

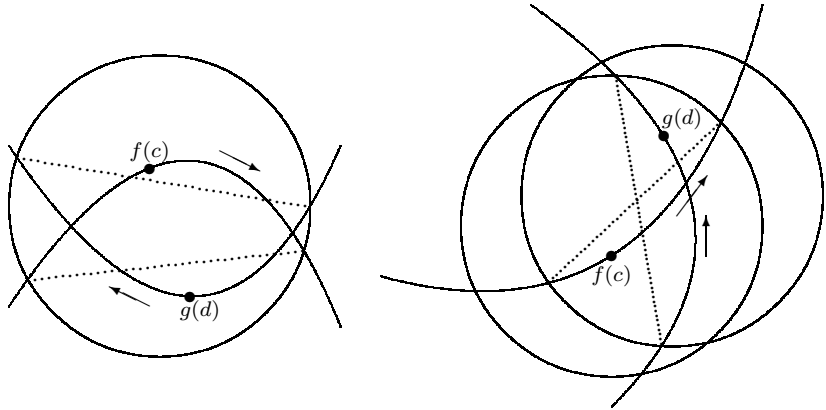
\begin{figure}[htbp]
\setlength{\unitlength}{3.8pt}
%\linethickness{0.7pt}
\begin{picture}(80,39)(-25,1)

\newsavebox{\ballu}
\savebox{\ballu}{
\qbezier(15,0)(15,6.21)(10.61,10.61)
\qbezier(10.61,10.61)(6.21,15)(0,15)
\qbezier(15,0)(15,-6.21)(10.61,-10.61)
\qbezier(10.61,-10.61)(6.21,-15)(0,-15)
\qbezier(-15,0)(-15,6.21)(-10.61,10.61)
\qbezier(-10.61,10.61)(-6.21,15)(0,15)
\qbezier(-15,0)(-15,-6.21)(-10.61,-10.61)
\qbezier(-10.61,-10.61)(-6.21,-15)(0,-15)
%\put(0,0){\circle*{.5}}
}

\put(15,20){\usebox{\ballu}}

\qbezier(0,10)(20,40)(33,8)
\put(14,23.7){\circle*{1}}
\put(12,24.7){\parbox{10ex}{ \footnotesize$f(c)$ }}
\put(21,25.5){\vector(2,-1){4}}
%\put(3.,11.3){\parbox{10ex}{\footnotesize $x_{Vc}$ }}
%\put(30.3,14){\parbox{10ex}{\footnotesize $y_{Vc}$ }}

\qbezier(0,26)(20,-4)(33,26)
\put(18,11){\circle*{1}}
\put(17,8.9){\parbox{10ex}{\footnotesize $g(d)$ }}
\put(14,10){\vector(-2,1){4}}
%\put(-3.5,23.5){\parbox{10ex}{\footnotesize $\overline y_{Vd}$ }}
%\put(30.3,18.3){\parbox{10ex}{\footnotesize $\overline x_{Vd}$ }}

\multiput(.8,24.8)(.7,-0.118){42}{\circle*{.2}}
\multiput(1.9,12.6)(.7,0.074){40}{\circle*{.2}}

%%%%%%%%%%%%%%%%%%%%%%%%%%%%%%%%%%%%%%%%%%%%%%%%%%%%%%%%%
\put(60,18){\usebox{\ballu}}
\put(66,21){\usebox{\ballu}}

\qbezier(52,40)(80,20)(60,0)
\qbezier(37,13)(66,5)(75,40)

%\put(58,31){\parbox{10ex}{\footnotesize $\overline y_{Vd}$ }}
%\put(61.8,7){\parbox{10ex}{\footnotesize $\overline x_{Vd}$ }}
%\put(51,10){\parbox{10ex}{\footnotesize $x_{Vc}$ }}
%\put(71.6,28){\parbox{10ex}{\footnotesize $y_{Vc}$ }}
\put(69.4,15){\vector(0,1){4}}

\put(65.2,27){\circle*{1}}
\put(65,28){\parbox{10ex}{\footnotesize $g(d)$ }}
\put(60,15){\circle*{1}}
\put(58,12.4){\parbox{10ex}{\footnotesize $f(c)$ }}
\put(66.5,19){\vector(3,4){3}}

\multiput(60.4,32.9)(.08,-0.467){58}{\circle*{.2}}
\multiput(53.75,12.42)(.4,.372){44}{\circle*{.2}}
\end{picture}
\caption{A ball-shaped gate not crossed by $g$ via $d$  and a lens-shaped gate crossed by $g$ via $d$.} \label{fig2}
\end{figure}

Figure~\ref{fig2} shows a ball-shaped gate $(V,c)$ which is not crossed by $g$ via $d$ since the dotted line segments  do not intersect (equivalently, surrounding the circle (boundary of $V$) the four endpoints of the curve segments do not alternate in $f$ and $g$).
Notice that nevertheless the curve segment $f(I_{Vc})$ of $f$ and the curve segment $g(I^g_{Vd})$ of $g$ intersect inside of $V$. Such intersection points will be irrelevant in our construction.

On the right Figure~\ref{fig2} shows a lens-shaped gate (intersection of two balls) $(V,c)$ which is crossed by the path $g$ via $d$ (the dotted line segments intersect).
As an example, $((0;1)^2,1/2)$  is a gate which is crossed by $g$ via the point $d=1/2$ (see Figure~\ref{fig1}). Here $I_{Vc}=I^g_{Vd}=(0;1)$

In the remainder of this section we prove the following central lemma which generalizes Step C of the proof of the computable intermediate value theorem in Section~\ref{seca}.
\begin{lemma}\label{l2} Suppose $\neg Q$ and let
$(V,c)$ be a gate crossed by $g$.
Then there is a gate $(V',c')$ crossed by $g$ such that $\overline {V'}\In V$ and   ${\rm diameter}(V')< {\rm diameter}(V)/2$.
\end{lemma}

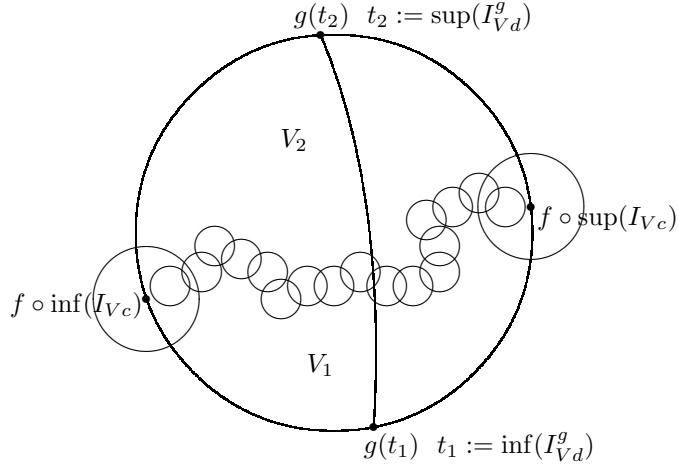
\begin{figure}[htbp]
\setlength{\unitlength}{5pt}
%\linethickness{0.7pt}
\begin{picture}(80,36)(-35,-2)

\newsavebox{\ballc}
\savebox{\ballc}{
\qbezier(15,0)(15,6.21)(10.61,10.61)
\qbezier(10.61,10.61)(6.21,15)(0,15)
\qbezier(15,0)(15,-6.21)(10.61,-10.61)
\qbezier(10.61,-10.61)(6.21,-15)(0,-15)
\qbezier(-15,0)(-15,6.21)(-10.61,10.61)
\qbezier(-10.61,10.61)(-6.21,15)(0,15)
\qbezier(-15,0)(-15,-6.21)(-10.61,-10.61)
\qbezier(-10.61,-10.61)(-6.21,-15)(0,-15)
%\put(0,0){\circle*{.5}}
}

\put(16,15){\usebox{\ballc}}

\qbezier(15,29.95)(20,18 )(19,0.27)

\put(1.8,10){\circle*{.5}}
\put(-3.5,9.5){\makebox(0,0)[cc]{$f\circ \inf(I_{Vc})$}}
\put(30.9,17){\circle*{.5}}
\put(36.7,16.1){\makebox(0,0)[cc]{$f\circ \sup(I_{Vc})$}}
\put(15,29.95){\circle*{.5}}
\put(22,31.4){\makebox(0,0)[cc]{$g(t_2)$\ \ $t_2:=\sup(I^g_{Vd})$}}
\put(19,0.27){\circle*{.5}}
\put(27,-1.2){\makebox(0,0)[cc]{$g(t_1)$\ \ $t_1:= \inf(I^g_{Vd})$}}

\put(1.8,10){\circle{10}}
\put(30.9,17){\circle{10}}

\put(3.6,11){\circle{3}}
\put(6,12){\circle{3}}
\put(7,14){\circle{3}}
\put(9,13){\circle{3}}
\put(11,12){\circle{3}}
\put(12,10){\circle{3}}
\put(14,11){\circle{3}}
\put(16,11){\circle{3}}
\put(18,12){\circle{3}}
\put(20,11){\circle{3}}
\put(22,11){\circle{3}}
\put(24,12){\circle{3}}
\put(24,14){\circle{3}}
\put(23,16){\circle{3}}
\put(25,17){\circle{3}}
\put(27,18){\circle{3}}
\put(29,17){\circle{3}}

\put(13,22){\makebox(0,0)[cc]{$V_2$}}
\put(15,5){\makebox(0,0)[cc]{$V_1$}}

\end{picture}
\caption{A ball-shaped gate $(V,c)$ with a chain of balls covering $f$ in $V$. } \label{fig4}
\end{figure}

We define a sequence $(V_i,a_i):=(B(f(a_i),r_i),a_i)$ of overlapping ball-shaped gates covering $f(\overline I_{Vc})$ such that also the intersections of consecutive balls (lenses) are gates  and show that at least one of these gates must be crossed by~$g$.
Figure~\ref{fig4} shows a ball-shaped gate $(V,c)$ with such a sequence of balls. The picture is topologically correct but not metrically. In general the sequence of balls may have numerous loops. The exact details are given below.

First we define  balls $V_0,V_1,\ldots,V_n,V_{n+1}$ and rational numbers $b_1,\ldots, b_{n-1}$.
Since the functions $f$  and  $g$ are continuous on the compact interval $[0;1]$, they are uniformly continuous. Therefore, there is an non-decreasing modulus function $\mod:\IQ\to\IQ$ such that for all $a,b\in [0,1]$,
 \begin{eqnarray}\label{f19}|f(a)-f(b)|<\gamma \ \ \mbox{and}\ \
 |g(a)-g(b)|<\gamma\ &
 \mbox{ if }\ & |a-b|<\mod(\gamma) \ \ \ \ (\mbox{for all } \  \gamma>0)\,.
 \end{eqnarray}
Since $f\circ \inf(I_{Vc}),f\circ \sup(I_{Vc})\not \in\range(g)$ (since $(V,c)$ is a gate) and $\range(g)$ is compact, there is a rational number  $\overline r$  such that
\begin{eqnarray}\label{f29}
B(f\circ \inf(I_{Vc}),2\overline r)\cap \range(g)=\emptyset \an B(f\circ \sup(I_{Vc}),2\overline r)\cap \range(g)=\emptyset
\end{eqnarray}
(see the first ball and the last ball in Figure~\ref{fig4}).
Since $(V,c)$ is crossed by $g$, \  $f\circ \inf(I_{Vc})\neq f\circ \sup(I_{Vc})$. Therefore  for avoiding pathological cases, we may assume
\begin{eqnarray}\label{f21}
B(f\circ \inf(I_{Vc}),2\overline r) \cap B(f\circ \sup(I_{Vc}),2\overline r)=\emptyset\,.
\end{eqnarray}
Since $f$ is continuous, there are rational numbers $a_l,a_r$ such that
\begin{eqnarray}\label{f17}
\inf(I_{Vc})<a_l<a_r<\sup(I_{Vc})\,,\ \ f[\inf(I_{Vc});a_l]\In B(f\circ \inf(I_{Vc}),\overline r) \an f[a_r;\sup(I_{Vc})]\In B(f\circ \sup(I_{Vc}),\overline r)\,.
\end{eqnarray}
Since $f[\inf(I_{Vc});a_l]$ and  $f[a_r; \sup(I_{Vc})]$ are compact and
$f[a_l,a_r]$ is a  compact subset of the open set~$V$, there is a rational number $s>0$ such that
\begin{eqnarray}
\label{f16}&& B(x,3s)\In  B(f\circ \inf(I_{Vc}),\overline r) \mbox{ for all } x\in f[\inf(I_{Vc});a_l]\,,\\
\label{f16a}&&B(x,3s)\In  B(f\circ \sup(I_{Vc}),\overline r) \mbox{ for all } x\in f[a_r; \sup(I_{Vc})] \an \\
\label{f30} &&\mbox{$B(x,3s)\In V$ \ for all $x\in f[a_l;a_r]$\,.}
\end{eqnarray}

For $i\geq 1$ let $\varepsilon _i:=2^{-i-1}\cdot s$. Then $\sum_{i\geq 1}\varepsilon_i=s/2$.
We choose $n\in\IN$, rational numbers $a_i,b_i,r_i,s_i$ and a real number $t_i $ inductively as follows.
\smallskip

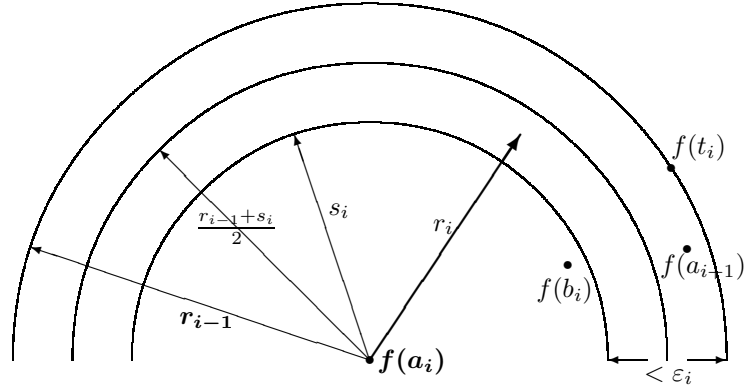
\begin{figure}[htbp]
\setlength{\unitlength}{.6pt}
%\linethickness{0.7pt}
\begin{picture}(500,250)(-440,-12)

\newsavebox{\balle}
\savebox{\balle}{
\qbezier(150,0)(150,62.1)(106.1,106.1)
\qbezier(106.1,106.1)(62.1,150)(0,150)
\qbezier(-150,0)(-150,62.1)(-106.1,106.1)
\qbezier(-106.1,106.1)(-62.1,150)(0,150)
\put(0,0){\circle*{5}}
}

\newsavebox{\ballf}
\savebox{\ballf}{
\qbezier(225,0)(225,93.15)(159.2,159.2)
\qbezier(159.2,159.2)(93.15,225)(0,225)
\qbezier(-225,0)(-225,93.15)(-159.2,159.2)
\qbezier(-159.2,159.2)(-93.15,225)(0,225)
}

\newsavebox{\ballg}
\savebox{\ballg}{
\qbezier(187.5,0)(187.5,77.63)(132.6,132.6)
\qbezier(132.6,132.6)(77.63,187.5)(0,187.5)
\qbezier(-187.5,0)(-187.5,77.63)(-132.6,132.6)
\qbezier(-132.6,132.6)(-77.63,187.5)(0,187.5)
}

\put(0,0){\usebox{\balle}}
\put(0,0){\usebox{\ballf}}
\put(0,0){\usebox{\ballg}}

\put(0,0){\vector(-3,1){214}}
\put(0,0){\vector(-1,1){132.7}}
\put(0,0){\vector(-1,3){47.5}}
\thicklines
\put(0,0){\vector(2,3){95}}
\thinlines

\put(175,0){\vector(-1,0){25}}
\put(200,0){\vector(1,0){25}}
\put(173,-14){\parbox{10ex}{ $<\varepsilon_i$ }}

\put(-120,20){\parbox{10ex}{\boldmath $r_{i-1}$ }}
\put(-26,90){\parbox{10ex}{ $s_i$ }}
\put(-110,80){\parbox{10ex}{ $\frac{r_{i-1}+s_i}{2}$ }}
\put(40,80){\parbox{10ex}{ $r_i$ }}

\put(5,-5){\parbox{10ex}{\boldmath $f(a_i)$ }}

\put(125,60){\circle*{5}}
\put(105,40){\parbox{10ex}{ $f(b_i)$ }}

\put(190,121){\circle*{5}}
\put(190,130){\parbox{10ex}{ $f(t_i)$ }}

\put(200,70){\circle*{5}}
\put(183,55){\parbox{10ex}{ $f(a_{i+1})$ }}

\end{picture}
\caption{The definition of $b_i$, $a_{i+1}$ and $r_i$ from $a_i$ and $r_{i-1}$. } \label{fig3}
\end{figure}

$r_0:= s$, $a_1:=a_l$.\\
Suppose $r_{i-1}$ and $a_i$ such that $r_{i-1}\leq s$ and $a_i\geq a_l$ have been chosen.\\
If $a_i\geq a_r$ then define $n:=i-1$ (End of construction.),
 else
choose  $t_i, b_i,s_i,a_{i+1}$ and $r_i$ such that
\begin{eqnarray}
\label{f48}&&t_i:=\sup\{b\in \IQ, a_i<b\mid f[a_i;b]\In B(f(a_i),r_{i-1})\}\,,\\
\label{f44}&&b_i\in\IQ\,, \ \ a_i<b_i \an t_i-\mod(\varepsilon_i)<b_i<t_i\,,\\
\label{f45}&&s_i<r_{i-1} \an f[a_i;b_i]\In B(f(a_i),s_i)\,,\\
\label{f46}&&  t_i-\mod((r_{i-1}-s_i)/2)<a_{i+1}<t_i\,,\\
\label{f47}&&s_i<r_i<(s_i+r_{i-1})/2 \an (B(f(a_i),r_i),a_i) \mbox{ is a gate}\,.
\end{eqnarray}
Figure~\ref{fig3} illustrates the construction. Notice that in general $\varepsilon_i$ is much smaller than shown in the figure since it converges to $0$ while $r_i >s/2$ for all $i$ (by (\ref{f52})).

We define
\begin{eqnarray}\label{f37}
 V_0:=B(f\circ \inf(I_{Vc}),\overline r)\,,& V_{n+1}:=B(f\circ \sup(I_{Vc}),\overline r)\,, &
V_i:=B(f(a_i),r_i) \ \ (1\leq i\leq n)\,.
\end{eqnarray}
These balls are shown in Figure~\ref{fig4}. We must show that the construction is sound, in particular that it ends after finitely many steps.
First we show that numbers $t_i, b_i,s_i,a_{i+1}$ and $r_i$ exist provided $a_i<a_r$.

{\boldmath $t_i$}: Since $a_l\leq a_i<a_r$,
$B(f(a_i),r_{i-1})\In B(f(a_i),s)\In V$ by (\ref{f30}). Therefore, $t_i$ exists.

{\boldmath $b_i$} exists.

{\boldmath $s_i$} exists since $b_i<t_i$.

{\boldmath $a_{i+1}$} exists since $s_i<r_{i-1}$

{\boldmath $r_i$}: Since $s_i<r_{i-1}$, $s_i<(s_i+r_{i-1})/2 $. By (\ref{f23})  $r_i$ exists such that $s_i<r_i<(s_i+r_{i-1})/2 $, \
 $f\circ \max(K_{a_ir_i})\not\in \range(g)$ and
  $f\circ \in(L_{a_ir_i})\not\in \range(g)$. By (\ref{f20}), $(B(f(a_i),r_i),a_i)$ is a gate.
\smallskip

Notice that (\ref{f47}) is the place where the assumption  $\neg Q$  is needed.
Below we will use the following equality which follows from  continuity of $f$:
\begin{eqnarray}\label{f8}
|f(a)-f(t)|=r & \mbox{if} & t=\sup\{b\in\IQ,\ a<b\mid f[a;b]\In
B(f(a),r)\} \ \ \ \mbox{and} \ \ \ B(f(a),r)\In (0;1)^2,.
\end{eqnarray}

\begin{lemma} \label{l1}
Suppose $i\geq 1$, $a_l\leq a_i<a_r$ and $s-\sum_{k=1}^{i-1}\varepsilon_k\leq r_{i-1}\leq s$. Then:
\begin{eqnarray}
\label{f51}&&r_i<r_{i-1}\,,\\
\label{f52}&&s-\sum_{k=1}^i\varepsilon_k<s_i\ \ (<r_i<r_{i-1}\leq s)\,,\\
\label{f53}&&f[a_i;b_i]\In B(f(a_i),r_i)\,,\\
\label{f54}&&f[b_i;a_{i+1}]\In B(f(a_{i+1}),s/2)\,,\\
\label{f55}&& r_i<|f(a_i)-f(a_{i+1})|<3s/2\,,\\
\label{f56}&&a_{i+1}\geq a_i+\mod(s/2)\,.
\end{eqnarray}
\end{lemma}
\proof $ $

{\bf(\ref{f51})} This follows from (\ref{f45}) and (\ref{f47}).

{\bf(\ref{f52})} By (\ref{f44}), $|f(t_i)-f(b_i)|<\varepsilon_i$ and
by (\ref{f45}), $|f(a_i)-f(b_i)|<s_i$. Therefore by (\ref{f8}) and (\ref{f47}),
$r_{i-1}=|f(a_i)-f(t_i)|\leq |f(a_i)-f(b_i)|+|f(b_i)-f(t_i)|
< s_i+\varepsilon_i$.

We obtain $s-\sum_{k=1}^i\varepsilon_k
=s-\sum_{k=1}^{i-1}\varepsilon_k-\varepsilon_i
\leq r_{i-1}-\varepsilon_i
<s_i<r_i<r_{i-1}<s$ \ (by (\ref{f47}) and (\ref{f51})).

{\bf(\ref{f53})} By (\ref{f45}) and (\ref{f47}), $f[a_i;b_i]\In B(f(a_i),s_i)\In B(f(a_i),r_i)$.

{\bf(\ref{f54})} Let $b_i\leq b\leq a_{i+1}$.\\

By (\ref{f44}) and (\ref{f46}),
 $t_i-\mod(\varepsilon_i)<b_i\leq b\leq a_{i+1}<t_i$, hence
$|a_{i+1}-b| <  \mod(\varepsilon_i)$, therefore by (\ref{f52}), $|f(b)-f(a_{i+1})|<\varepsilon_i <s/2$.
(\ref{f54}) follows immediately.

{\bf(\ref{f55})}
By (\ref{f52}--\ref{f54}),
$|f(a_i)-f(a_{i+1})|\leq |f(a_i)-f(b_i)| +|f(b_i)-f(a_{i+1})| <r_i+s/2
\leq 3s/2$.

  By (\ref{f46}), $|f(a_{i+1})-f(t_i)|<(r_{i-1} -s_i)/2$, hence by (\ref{f8}),
 $r_{i-1}=|f(a_i)-f(t_i)|\leq |f(a_i)-f(a_{i+1})| + |f(a_{i+1}) -f(t_i)|< |f(a_i)-f(a_{i+1})| + (r_{i-1} -s_i)/2$. Therefore,
 $ |f(a_i)-f(a_{i+1})| >r_{i-1}-(r_{i-1} -s_i)/2=(r_{i-1} +s_i)/2>r_i$
 by (\ref{f47}). \\
 This proves (\ref{f55}).

{\bf(\ref{f56})} Suppose $a_{i+1}-a_i<\mod(s/2)$. Then
$|f(a_{i+1})-f(a_i)|<s/2$. But by (\ref{f55}) and (\ref{f52}),
$|f(a_{i+1})-f(a_i)|>r_i\geq  s-\sum_{k=1}^i\varepsilon_k >s/2$.
\qq\\

Since the assumptions of Lemma~\ref{l1} are true for $i=1$ ($r_{i-1}=s$, $a_1=a_l$), the construction can be started with $i=1$.
If it has been performed for $i-1$ then  (\ref{f51})-(\ref{f56}) are true for $i-1$. By (\ref{f52}), $ r_{i-1}<s$  and $a_i\geq a_l$.
Then the construction can be continued with $i$, provided $a_i<a_r$.
By (\ref{f56}) there is a greatest number $i$ such that $a_i< a_r$. This number has been called $n$. Then $a_n<a_r\leq a_{n+1}$.

\begin{lemma}\label{l3}

For the balls $V_0,V_{n+1}$ and $V_i=B(f(a_i),r_i)$ ($\,1\leq i\leq n$) from  (\ref{f37}) (see Figure~\ref{fig4}),
 \begin{eqnarray}
\label{f22a}&&\overline  V_1\In V_0\,, \ \ \ \overline V_n\In V_{n+1}\,,\\
\label{f22b}&&\overline V_i\In B(f(a_i),s)\In B(f(a_i),3s)\In V \,.
 \end{eqnarray}
\end{lemma}

\proof $ \overline V_1=B(f(a_1),r_1)\in V_0$ by (\ref{f16}) since $a_1=a_l$ and $r_1<s$ by (\ref{f52}).
Since $B(f(a_{n+1}),3s)\In V_{n+1}$ (by (\ref{f16a})),
$|f(a_n)-f(a_{n+1})|<3s/2$  (by (\ref{f55})) and $r_n<s$, \
$\overline V_n\In B(f(a_n),s)\In B(f(a_{n+1}),3s)\In V_{n+1}$.
The third statement follows from   (\ref{f52}) and  (\ref{f30}). \qq
\medskip

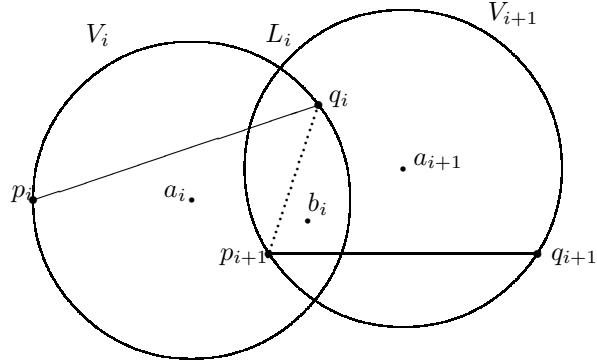
\begin{figure}[htbp]
\setlength{\unitlength}{4pt}
%\linethickness{0.7pt}
\begin{picture}(80,37)(-22,5)

\newsavebox{\ballk}
\savebox{\ballk}{
\qbezier(15,0)(15,6.21)(10.61,10.61)
\qbezier(10.61,10.61)(6.21,15)(0,15)
\qbezier(15,0)(15,-6.21)(10.61,-10.61)
\qbezier(10.61,-10.61)(6.21,-15)(0,-15)
\qbezier(-15,0)(-15,6.21)(-10.61,10.61)
\qbezier(-10.61,10.61)(-6.21,15)(0,15)
\qbezier(-15,0)(-15,-6.21)(-10.61,-10.61)
\qbezier(-10.61,-10.61)(-6.21,-15)(0,-15)
\put(0,0){\circle*{.5}}
}

\put(35,20){\usebox{\ballk}}
\put(32.4,20){\parbox{1ex}{ $a_i$ }}

\put(25,35){\parbox{1ex}{ $V_i$ }}
\put(63,37){\parbox{1ex}{ $V_{i+1}$ }}
\put(42,35){\parbox{1ex}{ $L_i$ }}

\put(55,23){\usebox{\ballk}}
\put(56,23){\parbox{1ex}{ $a_{i+1}$ }}

\put(20,20){\line(3,1){27}}
%\put(47,29){\line(-1,-3){4.69}}
\put(42.3,14.95){\line(3,0){25.4}}
\multiput(47,29)(-.2,-.6){25}{\circle*{.3}}

\put(20,20){\circle*{.7}}
\put(18,20){\parbox{1ex}{ $p_i$ }}
\put(47,29){\circle*{.7}}
\put(48,29){\parbox{1ex}{ $q_i$ }}
\put(42.3,14.95){\circle*{.7}}
\put(37.7,14){\parbox{1ex}{ $p_{i+1}$ }}
\put(67.7,14.95){\circle*{.7}}
\put(69,14){\parbox{1ex}{ $q_{i+1}$ }}

\put(46,18){\circle*{.5}}
\put(46,19){\parbox{1ex}{ $b_i$ }}

\end{picture}
\caption{The balls $V_i$ and $V_{i+1}$.\ The numbers $a_i,b_i,p_i,q_i$ are ``times" for $f$. } \label{fig8}
\end{figure}

For  $1\leq i<n$ let $L_i:=V_i\cap V_{i+1}$.
By (\ref{f51} -- \ref{f55}) for $1\leq i<n$,
\begin{eqnarray}\label{f36}
f(a_i)\not\in V_{i+1}\,,&f(a_{i+1})\not\in V_i\,, & f(b_i)\in V_i\cap V_{i+1}=L_i\,.
\end{eqnarray}
Therefore the positions of  $a_i,a_{i+1}$ and $b_i$ are drawn correctly in Figure~\ref{fig8}. By (\ref{f36})  $(L_i,b_i)$ is a passage. Define
$$(p_i,q_i):=I_{V_ia_i}\ \ \ \mbox{and}\ \ \ (p_{Li},q_{Li}):=I_{L_ib_i}\,.$$

The straight line segments in Figure~\ref{fig8} correspond to the dotted lines in Figure~\ref{fig2}.

\begin{lemma}\label{l8}
Suppose $1\leq i< n$. Then
\begin{enumerate}
\item \label{l8c}$(V_j,a_j)$ is a gate (for $1\leq j\leq 1$),
\item \label{l8a}$f(q_i)\in \partial V_i\cap V_{i+1}$ and
$f(p_{i+1})\in V_i\cap \partial V_{i+1}$\,,
\item \label{l8b}$q_i=q_{Li}$ and  $p_{i+1}=p_{Li}$\,.

 \item\label{l8d} $(L_i,b_i)$ is a gate.
\end{enumerate}
\end{lemma}

\proof $ $

\ref{l8c}: By (\ref{f47}).

\ref{l8a}: This follows from (\ref{f53})  and (\ref{f54})

\ref{l8b}: By (\ref{f53}) and  (\ref{f54}),
\begin{eqnarray*}
q_i&=&\sup\{b>a_i\mid f[a_i;b]\in V_i\}\\
&=&\sup\{b>b_i\mid f[b_i;b]\in V_i\}\\
&=&\sup\{b>b_i\mid f[b_i;b]\in V_i\cap V_{i+1}\}=q_{Li}\,,\\
p_{i+1}& = & \inf\{b<a_i\mid f[b;a_i]\in V_{i+1}\}\\
& = & \inf\{b<a_i\mid f[b;b_i]\in V_{i+1}\}\\
& = & \inf\{b<a_i\mid f[b;b_i]\in V_i\cap V_{i+1}\}=p_{Li}\,.
\end{eqnarray*}

\ref{l8d}: This follows from \ref{l8b}. and \ref{l8c}..
\qq\\

By Lemma~\ref{l8}.\ref{l8a} in Figure~\ref{fig8}  the points $f(p_{i+1})$ and $f(q_i)$ are positioned correctly, in particular by Lemma~\ref{l8}.\ref{l8a} they are different and not in
$\partial V_i\cap \partial V_{i+1}$ (the peaks of $L_i$).
Notice that $f(p_i)\in V_{i+1}$ is not excluded, but always $f(p_i)\neq f(q_i)$ by definition.
The line segments in Figure~\ref{fig8} correspond to the dotted line segments in Figure~\ref{fig2}.

\begin{figure}[htbp]
\setlength{\unitlength}{2.5pt}
%\linethickness{0.7pt}
\begin{picture}(100,65)(-52,-6)

\newsavebox{\ballv}
\savebox{\ballv}{
\qbezier(15,0)(15,6.21)(10.61,10.61)
\qbezier(10.61,10.61)(6.21,15)(0,15)
\qbezier(15,0)(15,-6.21)(10.61,-10.61)
\qbezier(10.61,-10.61)(6.21,-15)(0,-15)
\qbezier(-15,0)(-15,6.21)(-10.61,10.61)
\qbezier(-10.61,10.61)(-6.21,15)(0,15)
\qbezier(-15,0)(-15,-6.21)(-10.61,-10.61)
\qbezier(-10.61,-10.61)(-6.21,-15)(0,-15)
\put(0,0){\circle*{.9}}
}

\newsavebox{\xxx}
\savebox{\xxx}{
\qbezier(40.5,0)(40.5,16.77)(28.65,28.65)
\qbezier(28.65,28.65)(16.77,40.5)(0,40.5)
\put(0,0){\circle*{1}}
}

\qbezier(-3,0)(-8.5,20)(1,50)
\put(-9,-3){\parbox{20ex}{\footnotesize boundary of $V$}}

\put(27,0){\parbox{20ex}{\footnotesize boundary of $V_0$ }}

\put(-19,5){\parbox{10ex}{ $\inf(I_{Vc})$ }}
\put(-4.27,6){\usebox{\xxx}}
%%%%%%%%%%%%%%%%%%%%%%%%%%%%
\put(15,20){\usebox{\ballv}\circle*{.9}}
\put(12,21){\parbox{10ex}{ \footnotesize$a_1$ }}
\put(35,33){\usebox{\ballv}}
\put(58,30){\usebox{\ballv}}
\put(78,18){\usebox{\ballv}}
\put(98,20){\usebox{\ballv}}

\put(15,48){\parbox{10ex}{ $V_1$ }}
\put(33,52){\parbox{10ex}{ $V_2$ }}
\put(58,51){\parbox{10ex}{ $V_3$ }}
\put(79,45){\parbox{10ex}{ $V_4$ }}
\put(96,45){\parbox{10ex}{ $V_5$ }}
%%%%%%%%%%%%%%%%%%%%%%%%%5

\put(6,32){\circle*{1}}
\put(1,33){\makebox{ \footnotesize$p_1$ }}
\put(24,32){\circle*{1}}
\put(23,33.3){\makebox{ \footnotesize$q_1$ }}

\put(22.9,24){\circle*{1}}
\put(20.2,21){\makebox{ \footnotesize$p_2$ }}
\put(49.2,38){\circle*{1}}
\put(49,38.6){\makebox{ \footnotesize$q_2$ }}

\put(43.8,24.9){\circle*{1}}
\put(39,24){\makebox{ \footnotesize$p_3$ }}
\put(68.2,19){\circle*{1}}
\put(66.5,17){\makebox{ \footnotesize$q_3$ }}

\put(70.2,30.7){\circle*{1}}
\put(67,32.2){\makebox{ \footnotesize$p_4$ }}
\put(92.3,13){\circle*{1}}
\put(92.4,12.3){\makebox{ \footnotesize$q_4$ }}

\put(84.7,13){\circle*{1}}
\put(79.8,12){\makebox{ \footnotesize$p_5$ }}
\put(96,34.8){\circle*{1}}
\put(95,36.3){\makebox{ \footnotesize$q_5$ }}

\put(16.5,41.4){\makebox{ $L_1$ }}
\put(45,45){\makebox{ $L_2$ }}
\put(72,36){\makebox{ $L_3$ }}
\put(83,33){\makebox{ $L_4$ }}

\put(6,32){\line(1,0){18}}
\qbezier(24,32)(23.45,24)(22.9,24)

\qbezier(22.9,24)(36.05,31)(49.2,38)
\qbezier(49.2,38)(46.5,31.45)(43.8,24.9)
\qbezier(43.8,24.9)(56,21.95)(68.2,19)
\qbezier(68.2,19)(69.2,24.85)(70.2,30.7)
\qbezier(70.2,30.7)(81.25,21.85)(92.3,13)
\qbezier(92.3,13)(88.5,13)(84.7,13)
\qbezier(84.7,13)(90.35,23.9)(96,34.8)

\end{picture}
\caption{The balls $V_i=B(f(a_i),r_i)$, \  $i=1,\ldots,5$ inside of the set $V$ and the path $\rm Pf$.   } \label{fig6}
\end{figure}
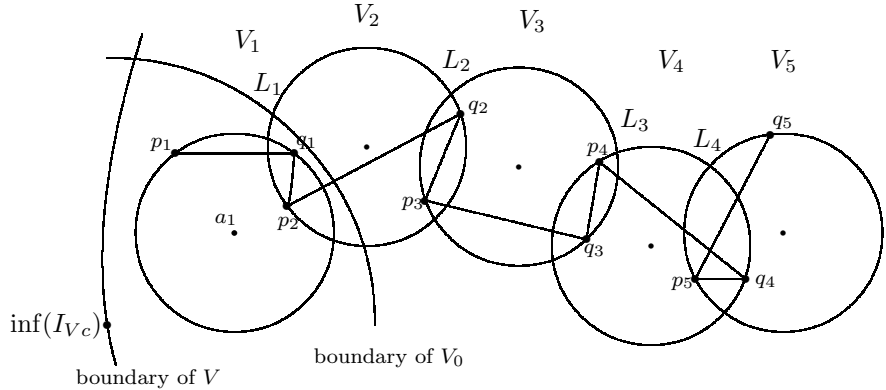

Figure~\ref{fig6} (as a detail of Figure~\ref{fig4})  shows an example of the balls  $V_1=B(f(a_1),r_1)$$\,\ldots,$ $V_5=B(f(a_5),r_5)$. The figure is topologically correct but not metrically, in particular since (\ref{f16}) and (\ref{f30}) are not satisfied.

For $1\leq i\leq n$  let $\gamma_i:=\overline{f(p_i)f(q_i)}$  be the straight line segment from $f(p_i)$ to $f(q_i)$ and for $1\leq i<n$ let
 $\delta_i:=\overline{f(q_i)f(p_{i+1})}$.
These line segments correspond to the dotted lines in Figure~\ref{fig2}.
The sequence $\gamma_1,\delta_1,\gamma_2,\ldots,\delta_{n-1},\gamma_n$
forms a polygon path $\rm Pf:=\bigcup \gamma_i$ running from  $p_1\in V\cap V_0$ to  $q_n \in V\cap  V_{n+1}$ (Figure~\ref{fig6}).\footnote{We have $p_i<a_i<p_{i+1}<b_i<q_i<a_{i+1}<q_{i+1}$  \ (where $p_{i+1}<b_i$ by (\ref{f54}) and $f(p_{i+1})\in V_i\setminus V_{i+1}$).
However, running along the straight line segments in Figure~\ref{fig8} the point $f(q_i)$ is visited before the point $f(p_{i+1})$. This is no contradiction.}

Now we complete the proof of Lemma~\ref{l2}. We have a gate $(V,c)$ crossed by $g$ via some $d$ as in Figure~\ref{fig4} ($f(c)$ and $g(d)$ are not marked in this figure) with more details in Figure~\ref{fig6}.
Let $(t_1,t_2):= I^g_{Vd}$. By (\ref{f29}),
 $$g(t_1),\ g(t_2)\in V':= \overline V\setminus(V_0\cup V_{n+1})\,.$$
\noindent The polygon ${\rm Pf}$ may have loops (not only loops inside a gate as shown in Figure~\ref{fig6} but much longer ones).
By cutting the loops we obtain a loop-free polygon ${\rm Pf'}\In {\rm Pf}$ running from $f(p_1)\in \partial V_1$ to $f(q_n)\in \partial V_n$.
Then $V'$ is the disjoint union of ${\rm Pf'}\cap V'$ and two connected sets $T_1,\,   T_2$ such that  the set $T_1 \cup T_2$ is not connected,
$g(t_1)\in T_1\cap \partial V$ and $g(t_2)\in T_2\cap \partial V$ since $(V,c)$ is crossed by $g$.
(The polygon ${\rm Pf'}$ cuts the set $V'$ into disjoint parts $T_l$ and $T_2$ with $g(t_1)\in T_1$ and $g(t_2)\in T_2$, see Figure~\ref{fig4}.)
\smallskip

Suppose ${\rm Pf}\cap g[t_1;t_2]=\emptyset$.   Then ${\rm Pf'}\cap g[t_1;t_2]=\emptyset$. Since $g(t_1)\in T_1$, $g[t_1;t_2]\in V'$ and
$g[t_1;t_2]$ is connected,  $g[t_1;t_2]\In T_1$. Therefore,
$$g[t_1;t_2]\In T_1\In T_1\cup {\rm Pf}\In T_1\cup \bigcup_{i=1}^n\overline  V_i\,.$$
We observe $g(t)$ for increasing $t\geq t_1$. Whenever $g$ enters a gate $(V_i,a_i)$ or $(L_i,b_i)$ it must leave it at the same side of the line $\gamma_i$ or $\delta_i$, respectively, without crossing it.
\medskip

Now we replace the assumption ${\rm Pf}\cap g[t_1;t_2]=\emptyset$ by the weaker assumption: none of the gates $(V_i,a_i)$ and $(L_i,b_i)$ is crossed by $g$. If $g$ enters $V_i$
(by Definition~\ref{d3} the endpoints of $\gamma_i$ are not in $\range(g)$) then it may cross the line $\gamma_i$ but must stay in $V_i$ and finally leave $V_i$ in the same side. The line segment $\gamma_i$ can be crossed only locally.
(correspondingly for $(L_i,b_i)$ and $\delta_i$), see the left example in Figure~\ref{fig2}.
Therefore, also with this weaker assumption
\begin{eqnarray}\label{f38} g(t_2)\in g[t_1;t_2]\In  T_1\cup \bigcup_{i=1}^n \overline  V_i\,.
\end{eqnarray}
Since by the definition of $T_1$ and $T_2$, $T_1\cap T_2=\emptyset$ and $g(t_2)\in T_2$, $g(t_2)\in T_1$ is impossible. By (\ref{f52}) and
(\ref{f30}),
$$d(\bigcup_{i=1}^n \overline  V_i,\partial V)\geq
d(\bigcup_{i=1}^n \overline B(f(a_i),s),\partial V)\geq 2s\,.$$
Since $g(t_2)\in\partial V$, $g(t_2)\in \bigcup_{i=1}^n \overline  V_i$ is impossible. Therefore (\ref{f38}) is false.
 Contradiction.\footnote
{Consider the path $\rm Pf$ as a fence with posts $p_1,q_1,p_2,q_2,\ldots, p_n,q_n$ and planks between consecutive posts.
Then it is not possible to walk from   $g(t_1)$ to $g(t_2)$  without crossing one of the planks. If we replace each plank by a belt which allows to cross the plank  only locally it is still impossible to walk from   $g(t_1)$ at time $g(t_2)$.}

\smallskip

We conclude that one of the gates $(V_i,a_i)$ and $(L_i,b_i)$ must be crossed by $g$. This ends the proof of Lemma~\ref{l2}.

\section{The second alternative}\label{sece}
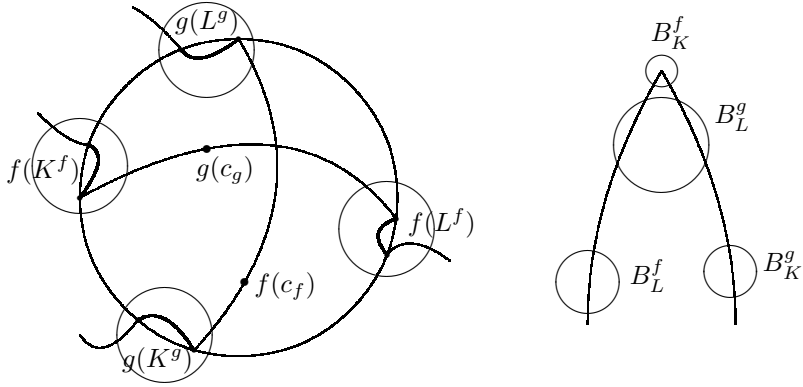
\begin{figure}[htbp]
\setlength{\unitlength}{4pt}
%\linethickness{0.7pt}
\begin{picture}(30,39)(-10,-18)

\newsavebox{\ballt}
\savebox{\ballt}{
\qbezier(15,0)(15,6.21)(10.61,10.61)
\qbezier(10.61,10.61)(6.21,15)(0,15)
\qbezier(15,0)(15,-6.21)(10.61,-10.61)
\qbezier(10.61,-10.61)(6.21,-15)(0,-15)
\qbezier(-15,0)(-15,6.21)(-10.61,10.61)
\qbezier(-10.61,10.61)(-6.21,15)(0,15)
\qbezier(-15,0)(-15,-6.21)(-10.61,-10.61)
\qbezier(-10.61,-10.61)(-6.21,-15)(0,-15)
}
\put(40,0){\usebox{\ballt}}
\put(25,3){\circle{9}}
\put(54,-3){\circle{9}}
\put(37,14){\circle{9}}
\put(33,-13){\circle{9}}

\put(25,0){\circle*{.5}}
\put(25.8,5){\circle*{.5}}

\put(54.05,-5.4){\circle*{.5}}
\put(54.9,-2){\circle*{.5}}

\qbezier(25,0)(45,11)(54.9,-2)
\thicklines
\qbezier(25,0)(28,4)(25.8,5)
\qbezier(54.9,-2)(52,-3)(54.05,-5.4)
\thinlines
\qbezier(25.8,5)(24,5)(21,8)
\qbezier(54.05,-5.4)(56,-3)(60,-6)

\put(40,15){\circle*{.5}}
\put(34.6,14){\circle*{.5}}

\put(30.5,-11.6){\circle*{.5}}
\put(35.8,-14.45){\circle*{.5}}

\qbezier(40,15)(49,0)(35.8,-14.45)
\thicklines
\qbezier(40,15)(36,12)(34.6,14)
\qbezier(35.8,-14.45)(33,-10)(30.5,-11.6)
\thinlines
\qbezier(34.6,14)(32,17)(30,18)
\qbezier(30.5,-11.6)(27,-16)(25,-13)

\put(40.57,-8){\circle*{.8}}
\put(41.5,-9){\parbox{10ex}{$\small f(c_f)$ }}

\put(37,4.6){\circle*{.8}}
\put(36,2){\parbox{10ex}{$\small g(c_g)$ }}

\put(18,2){\parbox{10ex}{$\small f(K^f)$ }}
\put(56,-3){\parbox{10ex}{$\small f(L^f)$ }}
\put(29,-16){\parbox{10ex}{$\small g(K^g)$ }}
\put(34,16){\parbox{10ex}{$\small g(L^g)$ }}

\qbezier(73,-12)(73,0)(80,12)
\qbezier(87,-12)(87,0)(80,12)

\put(80,12){\circle{3}}
\put(80,5){\circle{9}}
\put(73,-8){\circle{6}}
\put(86.5,-7){\circle{5}}

\put(79,15){\parbox{10ex}{$\small B^f_K$ }}
\put(85,7){\parbox{10ex}{$\small B^g_L$ }}

\put(77,-8){\parbox{10ex}{$\small B^f_L$ }}
\put(89.6,-7){\parbox{10ex}{$\small B^g_K$ }}

\end{picture}\label{fig7}
\caption{A proper ball-shaped proper crossing and a ball $B^g_L$ violating (\ref{f33})}
\end{figure}

In this section we generalize Step~{\bf D} of the proof of the computable intermediate value theorem in Section~\ref{seca}. We generalize the sets $K$ and $L$ from (\ref{f7}), (\ref{f15}) and Figure~\ref{fig5} to gates.
If $(V,c)$ is a gate crossed by $g$ via $d$ we call $(V,c,d)$ a crossing.
\begin{definition}\label{d4}$ $
\begin{enumerate}
\item \label{d4a}A {\em meeting} is a triple $(V,c_f,c_g)$, such that $c_f,c_g\in\IQ$ and $f(c_f),g(c_g)\in V$ where $V$ is either a ball $B(f(a),r)$ (a ball-shaped meeting) such that $a,r\in\IQ$ and  $\overline{B(f(a),r)}\In (0,1)^2$ or $V$ is the strict intersection\footnote{We call the intersection of $B_1$ and $B_2$ strict if $B_1\not\In B_2$, \ $B_2\not\In B_1$ and $B_1\cap B_2\neq\emptyset$.} of two such balls, $V=B_1\cap B_2$ where $B_1=B(f(a_1),r_1)$ and $B_2= B(f(a_2),r_2)$ (a lens-shaped meeting).
    For a lens-shaped meeting let $\{z_1,z_2\}:= \partial B_1\cap \partial B_2$ (the peaks of $V$).

 \item \label{d4b}For a meeting  $(V,c_f,c_g)$ define
\begin{eqnarray*}
 K^f&:=&[p'_f,p_f]\ \ \ \mbox{and}\\
L^f&:=&[q_f,q'_f]\ \ \ \mbox{where}\\
(p_f\,;q_f)& :=&\mbox{the longest open interval $I$ such that }\ c_f\in I \ \mbox{ and }\ f(I)\In V\,,\ \ \mbox{and}  \\
(p'_f;q'_f)&:=&  \mbox{the longest open interval $I$ such that }\ c_f\in I \ \mbox{ and }\ f(\overline I)\In \overline V\,.
\end{eqnarray*}
The numbers and intervals $p_g,q_g,p'_g,q'_g,K^g$ and $L^g$ are defined accordingly.
\item \label{d4c}A meeting  $(V,c_f,c_g)$  is a crossing if
the four points $f(p_f),f(q_f), g(p_g)$ and $g(q_g)$
 are pairwise different and alternate in $f$ and $g$ on the boundary of $\,V$
(equivalently, if the straight line segments from $f(p_f)$ to  $f(q_f)$ and from $g(p_g)$ to  $g(q_g)$ intersect in $V$).

\item \label{d4d}A meeting $(V,c_f,c_g)$ is a proper crossing, if there are rational balls (i.e. balls $B(c,r)$ with rational $c$ and~$r$) $B^f_K,B^f_L,B^g_K,B^g_L$ such that
\begin{eqnarray}
\label{f58}&&\mbox{the four balls are pairwise disjoint},\\
\label{f49} &&\mbox{the four balls alternate in $f$ and $g$ on the boundary of $V$\,}\,,\\
\label{f39}
&&f(K^f)\In B^f_K\,, \ \  f(L^f)\In B^f_L\,,\ \
f(K^g)\In B^g_K\,, \ \  f(L^g)\In B^g_L \,.
\end{eqnarray}
and if $\,V=B_1\cap B_2$ is lens-shaped then for each of the four balls $B$,
\begin{eqnarray}\label{f33}
 \overline B\In B_1\ \ \mbox{or}\ \ \overline B\In B_2,\ \ \mbox{or}\ \ z_1\in B\ \ \mbox{or}\ \ z_2\in B\,.
\end{eqnarray}
\end{enumerate}
\end{definition}
Figure~\ref{fig7} shows a ball-shaped proper crossing. The right side shows a ball $B^g_L$ violating~(\ref{f33}).

\begin{lemma}\label{l7}
Every proper crossing is a crossing.
\end{lemma}

\proof By the definitions $f(p_f)\in f(K^f)\cap\partial V$,
$f(q_f)\in f(L^f)\cap\partial V$,
$g(p_g)\in g(K^g)\cap\partial V$ and
$g(q_g)\in g(L^g)\cap\partial V$ (see Figure~\ref{fig5}).
Since the four sets are pairwise disjoint the four points are pairwise different.
By (\ref{f49})
the four points $f(p_f),f(q_f), g(p_g)$ and $g(q_g)$
alternate in $f$ and $g$ on the boundary .
 of~$\,V$
 \qq

\medskip
\noindent In our new terminology Lemma~\ref{l2} can be expressed as follows:

\begin{lemma}\label{l6} For every crossing $(V,c_f,c_g)$ there is a crossing $(V',c'_f,c'_g)$ such that $\overline{V'}\In V$ and
${\rm diameter}(V')< {\rm diameter}(V)/2$.
\end{lemma}

We will compute an intersection point of $f$ and $g$ as the limit of a nested decreasing sequence of crossings. Unfortunately the set of crossings is not c.e. since the endpoints of the sets $K$ and $L$ are not
computable so that the condition from Definition~\ref{d4}.\ref{d4c} is not c.e.\footnote{For example, for a crossing $(V,c_f,c_g)$  possibly
$f(K^f)= f(L^f) = f(K^g)= f(L^g)=\partial V$.}
We solve the problem by Lemma~\ref{l9} and Lemma~\ref{l10} below.

\begin{lemma}\label{l9}
For every  crossing $(V,c_f,c_g)$ there is a proper crossing $(W,c_f,c_g)$
such that $W\In V$.
\end{lemma}

\proof
We consider the case of a lens-shaped crossing.  The proof for ball-shaped crossings is similar but simpler. Consider a crossing $(V,c_f,c_g)$ where $V=B_1\cap B_2$. For $\delta>0$ let
$V_\delta:=B_{1\delta}\cap B_{2\delta}$ where $B_{1\delta}:= B(f(a_1),r_1-\delta)$ and $B_{2\delta}:= B(f(a_2),r_2-\delta)$.

There is some $\delta_0>0$ such that $\{f(c_f),g(c_g)\}\In V_\delta$ if $0\leq \delta < \delta_0$. Therefore, for $0\leq \delta< \delta_0$
we have intervals and points $K^f_\delta,\ldots, q'_{g\delta}$ according to Definition~\ref{d4}.\ref{d4b} (see Figure~\ref{fig5}).

\bigskip
By Definition~\ref{d4}.\ref{d4c}
the four points  $f(p_f),f(q_f), g(p_g)$ and $g(q_g)$ are pairwise different.
Therefore, there is some $e\in \IQ$, $e>0$,   such that the four balls
$B(f(p_f),2e), B(f(q_f),2e),B(g(p_g),2e)$ and $ B(g(q_g),2e)$
are pairwise disjoint and  for each of these balls $B$ with radius $2e$,\ \ %
 \begin{eqnarray}\label{f41}\overline B\In B_1\ \ \mbox{or}\ \ \overline B\In B_2\ \ \mbox{or} \ \ z_1= {\rm center}(B) \ \ \mbox{or}\ \ z_1= {\rm center}(B)\,.
 \end{eqnarray}
First we consider the point $q_f$.
Since $f$ is continuous there is some $s\in\IQ$, $c_f<s<q_f$, such that
$f[s;q_f]\In B(f(q_f),e)$. There is some $0<\delta_1<\min(\delta_0,e/2)$ such that $f[c_f;s]\In V_{\delta_1}$.
For every $0<\beta<\delta_1$ we obtain  $s<L^f_\beta<q_f$
(where $L^f_\beta:=$ the longest interval $I$ such that $c_f\in I$ and $f(I)\In V_\beta$), hence
\begin{eqnarray}\label{f57}
f(L^f_\beta)\In f[s;q_f]\In B(f(q_f),e)& \mbox{if} & 0<\beta<\delta_1,.
\end{eqnarray}
Correspondingly for the other numbers $p_f,p_g,q_g$ there are numbers $\delta_2,\delta_3$ and $\delta_4$.
Choose some $\gamma\in\IQ$ such that
\begin{eqnarray}\label{f50}0<\gamma< \min\{\delta_1,\delta_2,\delta_3,\delta_4\},\  \
|z_{1\gamma}-z_1|< e/2,\ \ |z_{2\gamma}-z_2|< e/2\,.
\end{eqnarray}

We show that $(V_\gamma,c_f,c_g)$ is a proper crossing.
Since $\gamma<\delta_1<\delta_0$, \ $\{f(c_f),g(c_g)\}\In V_\gamma$.
Therefore, $(V_\gamma,c_f,c_g)$ is a meeting.
For $(V_\gamma,c_f,c_g)$  let $z_{1\gamma},z_{2\gamma},p_{f\gamma},\ldots, K^f_\gamma,\ldots $ be the numbers and intervals introduced in
Definition~\ref{d4}.
We must verify Definition~\ref{d4}.\ref{d4d} for $(V_\gamma,c_f,c_g)$.
 By (\ref{f57}),
\begin{eqnarray}\label{f40}   f(K^f_\gamma)\In B(f(p_f),e),\ \  f(L^f_\gamma)\In B(f(q_f),e),\ \
f(K^g_\gamma)\In B(f(p_g),e),\ \  f(L^g_\gamma)\In B(f(q_g),e)\,,
\end{eqnarray}
where the balls are not necessarily rational.

First we  consider $L^f_\gamma$.
Since $f(L^f_\gamma)$ is a compact subset of $ B(f(q_f),e)$ there is some $c\in\IQ^2$ such that
\begin{eqnarray}\label{f42} |f(q_f)-c|<e/2& \mbox{and} & f(L^f_\gamma)\In B(c,e)\In B(f(q_f),3e/2)\,.
\end{eqnarray}
Define $B^f_{L\gamma}:= B(c,e)$. Then $f(L^f_\gamma)\In B^f_{L\gamma}$. Therefore (\ref{f39}) is true for $L^f_\gamma$.
We apply (\ref{f41}).

Suppose $\overline{B(f(q_f),2e)}\In B_1$.
Since $|f(q_f)-c|<e/2$ and $\gamma< e/2$, \ $B^f_{L\gamma}=B(c,e)\In
B(f(q_f),3e/2)\In B_{1\gamma}$. Therefore (\ref{f33}) is true for
the case $V_\gamma$. If   $\overline{B(f(q_f),2e)}\In B_2$ then
(\ref{f33}) is true for the case $V_\gamma$ accordingly.

Suppose in (\ref{f41}) $z_1={\rm center}(B(f(q_f),2e))$, that is, $f(q_f)=z_1$. By  (\ref{f50}) and (\ref{f42}),
$|z_{1\gamma}-c|\leq |z_{1\gamma}-z_1| + | z_1-c| <e/2+e/2=e\,,$
hence $z_{1\gamma}\in B(c,e)= B^f_{L\gamma}$.  Therefore (\ref{f33}) is true for the case $V_\gamma$. If $z_2={\rm center}(B(f(q_f),2e))$
then (\ref{f33}) is true for the case $V_\gamma$ accordingly.

All of this can be proved for $K^f_\gamma$, $K^g_\gamma$ and $L^g_\gamma $ (see (\ref{f40})) accordingly. Therefore (\ref{f33}) has been proved for $(V_\gamma,c_f,c_g)$.

Since $(V,c_f,c_g)$ is a crossing the points  $f(p_f),f(q_f), g(p_g)$ and $g(q_g)$ are pairwise different and alternate in $f$ and $g$ on the boundary of $\,V$. Let $u>0$ be be a lower bound of the mutual distances of these four points. We may choose $e< u/100$ and $\gamma< u/100$.
Then obviously the four balls
$B^f_{K\gamma},B^f_{L\gamma},B^g_{K\gamma}$ and $B^g_{L\gamma}$
alternate in $f$ and $g$ on the boundary of $V_\gamma$. We omit a detailed
verification. This proves (\ref{f49}) for $V_\gamma$.

Therefore, $(V_\gamma,c_f,c_g)$ is a proper crossing
\qq

\begin{lemma}\label{l10}
The set of ball-shaped proper crossings and the set of lens-shaped proper crossings are c.e.
\end{lemma}

More precisely, the set of tuples $(a,r,c_f,c_g)\in\IQ^4$ such that
$(B(f(a),r),c_f,c_g)$ is a proper crossing is c.e., and the set of tuples $(a_1, r_1, a_2,r_2,c_f,c_g)\in\IQ^6$ such that
$(B(f(a_1),r_1)\cap B(f(a_2),r_2),c_f,c_g)$ is a proper crossing is c.e.
\medskip

\proof We consider the case of lens-shaped crossings. The proof for ball-shaped crossings is similar. We use Definition~\ref{d4}.

For a tuple $(a_1, r_1, a_2,r_2,c_f,c_g)\in\IQ^6$,\ \,
$(B(f(a_1),r_1)\cap B(f(a_2),r_2),c_f,c_g)$ is a proper crossing if
there are rational balls $B^f_K,\,B^f_L,\,B^g_K$ and $B^g_L$ such that
$(V,c_f,c_g)$ is a meeting and (\ref{f58}) - (\ref{f33}) are true.

The predicate ($(V,c_f,c_g)$ is a meeting) is c.e. The properties (\ref{f58}), (\ref{f49}) and (\ref{f33}) are decidable.
Since $K^f\In (d_1;d_2)$ (for rational $d_1, d_2$) is c.e. (c.f. Lemma~\ref{l4}), $K^f$ can  computed (as a compact set). Since $f$ is computable, $f(K^f)$ can be computed (as a compact set). Therefore
$f(K^f)\In B^f_K$ is c.e. Since  the other three predicates from
(\ref{f39}) are c.e. accordingly, \  (\ref{f39}) is c.e.

Therefore,  the set of tuples $(a_1, r_1, a_2,r_2,c_f,c_g)\in\IQ^6$ such that
$(B(f(a_1),r_1)\cap B(f(a_2),r_2),c_f,c_g)$ is a proper crossing is c.e.
\qq\\

We can now prove the generalization of Step D of the proof of the computable intermediate value theorem in Section ~\ref{seca}.

\begin{theorem} \label{t5} If\, $\neg Q$ then there is a computable point $x\in \range(f)\cap\range(g)$.
\end{theorem}

\proof  We assume $\neg Q$ (\ref{f43}).
By Definition~\ref{d3} the unit square, more precisely $((0;1)^2,1/2)$ is a gate which is crossed by $g$ via $1/2$. By Lemma~\ref{l2} there is a gate $(V',c')$ crossed by $g$ via some $d$ such that ${\rm diameter}(V')<1$.
According to Definition~\ref{d4},$(V',c',d)$ is a crossing. which by Lemma~\ref{l9} has a smaller proper subcrossing.
By Lemma~\ref{l6} and Lemma~\ref{l9} every proper crossing has a has a proper subcrossing  with much smaller diameter.
By Lemma~\ref{l10} we can compute a sequence $((V_i,c_{fi},c_{gi}))_{i\in\IN}$ of proper crossings such that $\overline V_{i+1}\In V_i$ and ${\rm diameter}(V_i)<2^{-i}$.
The limit $x\in\IR^2$ such that $\{x\}=\bigcap _iV_i$ is a computable point.

Every neighborhood of $x$ contains $f(c_{fi})\in\range(f)$ for some $i \in \IN$. Since $[0;1]$ is compact, $\range(f)$ is compact, hence complete. Therefore, $x\in\range(f)$. Accordingly, $x \in \range(g)$.
\qq\\

The main theorem~\ref{t6} follows from Theorems~\ref{t4} and~\ref{t5}.

\section{A special case and additional remarks}\label{secf}
The intersection problem becomes easier in special cases.
\begin{theorem}\label{t7} \cite{IP17}
Suppose, in Theorem~\ref{t6} additionally the function $f$ is assumed to be injective. Then there are  computable real numbers $\alpha$ and $\beta$ such that $f(\alpha)= g(\beta)=x$.
\end{theorem}
\proof  Instead of the predicates $Q$ (\ref{f43}) and $Q_{iv}$ (\ref{f61}) we use \
$ Q_1:\iff (g(I)\In\range(f)\ \  \ \mbox{for some open interval}\ \ I)\,.$

If $Q_1$ then $g(\beta)\in \range(f)$ for some $\beta\in\IQ$ (which is computable).

Suppose $\neg Q_1$. Since $f$ is injective $\range(f)\cup \gamma_1\cup \gamma_2$ is a Jordan curve, where $\gamma_1$ and $\gamma_2$ are the straight line segments from $(0,0)$ to $(0,1)$ and from $(0,1)$ to $(1,1)$, respectively. Let $U$ be the interior of this curve and $V$ be its exterior.

Let $0<a<b<1\,$ be rational numbers such that $g(a)\in U$  and $g(b)\in V$ and let $\varepsilon >0$.
Since $g[a;b]$  is connected and $U\cup V$ is not connected, $g(t)\in\range(f)$ for some $a<t<b$. Let $t_0:=\sup\{t>a\mid g(t)\not\in V\}$. By the definition of $t_0$ there is some $b_1\in\IQ$ such that $t_0<b_1<t_0+\varepsilon <b$ and $g(b_1)\in V$. Since $\neg Q$, there is some $a_1\in\IQ $ such that $a<t_0-\varepsilon<a_1<t_0$ and $g(a_1)\in U$.
Therefore, there are rational numbers $a_1,b_1 $ such that $a<a_1<b_1<b$, $g(a_1)\in U$, $g(b_1)\in V$ and $b_1-a_1<2\varepsilon$.

The open sets $U$ and $V$ are computable (see Footnote~2). Therefore the sets $\{a\in\IQ\mid g(a)\in U\}$ and $\{a\in\IQ\mid g(a)\in V\}$ are c.e.
We can find rational numbers $a_0<b_0$ such that $g(a_0)\in U$ and $g(b_0)\in V$.
We can compute a nested sequence $(a_i;b_i)_i$ of rational intervals converging to some $\beta$ (which is computable) and such that that $g(a_i)\in U$ and $g(b_i)\in V$. Therefore, $g(\beta)\not\in V$ and $g(\beta)\not\in U$, hence $g(\beta)\in\range(f)$.

 Since $f$ is injective, in both cases from $\beta$ we can compute the real number $\alpha$ such that $f(\alpha)=g(\beta)$.
\qq\\

Theorem~\ref{t4} and Corollary~\ref{cor1} provide additional information to Theorem~\ref{t6}:  If $Q$ then there are  computable real numbers $\alpha$ and $\beta$ such that $f(\alpha)\in \range(g)$ and $g(\beta)\in\range (f)$. But the proof does not show that there are  computable real numbers $\alpha$ and $\beta$ such that $f(\alpha)=g(\beta)$. Accordingly, the proof of Theorem~\ref{t5} does not show that there are computable numbers $\alpha,\beta$ such that $f(\alpha)=x$ or $g(\beta)=x$.

In general, from the fact that $x\in\range (g) $ is computable we cannot conclude that $g(t)=x$ for some computable number $t$.
In the following example $g(t)=x$ for many real numbers $t$ but for no computable one.

\begin{example}\rm
\label{ex1} Let $x:= (1/2,1/2)$. Let $h_0$ be the function from Footnote~1\,. From $h_0$ we can define a computable function $h_1$ such that $h_1(1/3)=-1/2=h_1(2/3)$ and  $h_0$ and $h_1$ have the same zeroes.
Define

$$ g(t):=\left\{
\begin{array}{lll}
(3/2t,0)& \mbox{if} & 0\leq t\leq 1/3\,,\\
(1/2,h_1(t)+1/2) & \mbox{if} & 1/3\leq t\leq 2/3\,,\\
((3t-1)/2,0) & \mbox{if} & 2/3\leq t\leq 1\,.
\end{array}
\right .$$
Then $g(t)=x$ for many non-computable numbers $t$ but for no computable one.\qq
\end{example}

If $f$ and $g$ are are continuous but not necessarily computable our constructions are still computable in $f$ and $g$ (w.r.t. the canonical representation of the continuous functions $h:[0;1]\to [0;1]^2$ \cite{Wei00,WG09}). The predicate $Q$ (\ref{f16}) however, is not continuous (hence not computable) in $f,g$) and furthermore, in the proof of Theorem~\ref{t4}, the decision required in (\ref{f60}) is is not continuous (hence not computable) in $f,g$.  Possibly, the predicate $Q$ is not optimal. Also for finding $\alpha$ and $\beta$ such that $f(\alpha)=g(\beta)$ possibly some other predicate must be used

The degrees of unsolvability of
the operator $(f,g)\mmto x$ and the operator $(f,g)\mmto (\alpha,\beta)$ such that $f(\alpha)=g(\beta)$
have not yet been located in the Weihrauch lattice \cite{BGP17u}.
 Let ${\rm IP1}$ be the operator which finds a point of intersection from continuous $f$ and $g$ such that $f$ is injective (the uniform problem of Theorem ~\ref{t7}). Then
$\rm CC_{[0:1]}  \equiv_{sW} IVP \equiv_{sW}IP1$
\cite[Theorem~7.34]{BGP17u}.

\bibliographystyle{plain}

%\bibliography{cca_2017-06-20,meinebib_2018-02-19}
\end{document}